\documentclass[10pt,a4paper,twoside,openany]{report}
\usepackage[cp1251]{inputenc} 
\usepackage{amsfonts,amsmath,amssymb,amscd}
\usepackage{bezier,graphics,graphicx}
\usepackage{epsfig,layout,latexsym}
\usepackage{theorem}
\usepackage{cite}
\usepackage{indentfirst}
\usepackage{array}

\AtBeginDocument{} \AtBeginDocument{}

\theorembodyfont{\sl}
\newtheorem{theorem}{Theorem}[section]

\newtheorem{lemma}[theorem]{Lemma}
\newtheorem{corollary}[theorem]{Corollary}
\newtheorem{remark}[theorem]{Remark}

\numberwithin{equation}{section} \numberwithin{theorem}{section}
\newcommand{\card}{\text{card}\medskip}

\makeatletter
\renewenvironment{thebibliography}[1]
{\section*{\centerline{\rm\textsc{Bibliography}}}%
	\@mkboth{\MakeUppercase\refname}{\MakeUppercase\refname}%
	\list{\@biblabel{\@arabic\c@enumiv}}%
	{\settowidth\labelwidth{\@biblabel{#1}}%
		\leftmargin\labelwidth
		\advance\leftmargin\labelsep
		\@openbib@code
		\usecounter{enumiv}%
		\let\p@enumiv\@empty
		\renewcommand\theenumiv{\@arabic\c@enumiv}}%
	\sloppy
	\clubpenalty4000
	\@clubpenalty \clubpenalty
	\widowpenalty4000%
	\sfcode`\.\@m
	\setlength{\itemsep}{-0.1cm}}
{\def\@noitemerr
	{\@latex@warning{Empty 'thebibliography' environment}}%
	\endlist}
\renewcommand{\@biblabel}[1]{#1.}

\makeatother \hfuzz=0.5pt \tolerance=500

\begin{document}



\begin{center}
	\textbf{An optimal method for high order mixed derivatives of bivariate functions}

	\setcounter{footnote}{2} 
\end{center}


\vspace*{5mm}
\centerline{\textsc { Y.V. Semenova $\!\!{}^{\dag}$, S.G. Solodky  $\!\!{}^{\dag,\ddag}$} }


%
\vspace*{5mm}
\centerline{$\!\!{}^{\dag}\!\!$ Institute of Mathematics, National Academy of Sciences of Ukraine, Kyiv}
\centerline{$\!\!{}^{\ddag}\!\!$ University of Giessen, Department of Mathematics, Giessen, Germany}


	
\begin{small}
	\begin{quote}
	The problem of optimal recovering high-order mixed derivatives of bivariate functions with finite
	smoothness is studied. Based on the truncation method, an algorithm for numerical differentiation is
	constructed, which is order-optimal both in the sense of accuracy and in terms of the amount of involved Galerkin
	information. Numerical examples are provided to illustrate the fact that our approach can be implemented
	successfully.
	\end{quote}
\end{small}



\section{Description of the problem}\label{DescrProb}

Currently, many research activities on the problem of stable numerical differentiation have been taking place due to
the importance of this tool in such areas of science and technology as finance,
mathematical physics, image processing, analytical chemistry, viscous elastic mechanics, reliability analysis, pattern
recognition, and many others. Among these investigations, we highlight \cite{Dolgopolova&Ivanov_USSR_Comput_Math_Math_Phys_1966_Eng},
which is the first publication on numerical differentiation in terms of the theory of ill-posed problems.
Further research \cite{Dolgopolova&Ivanov_USSR_Comput_Math_Math_Phys_1966_Eng} has been continued in numerous
publications on numerical differentiation (see, for example, \cite{Ramm_1968_No11}, \cite{VasinVV_1969_V7_N2},
\cite{Groetsch_1992_V74_N2}, \cite{Hanke&Scherzer_2001_V108_N6}, \cite{Ahn&Choi&Ramm_2006}, \cite{Lu&Naum&Per},
\cite{Nakamura&Wang&Wang_2008}, \cite{Zhao_2010}, \cite{Zhao&Meng&Zhao&You&Xie_2016}, \cite{Meng&Zhaoa&Mei&Zhou_2020},
\cite{SSS_CMAM}),  covering different classes of the functions and the types of proposed methods. Despite the abundance
of works on this topic, the problem for recovery of high-order derivatives was considered only in a few publications,
among which we note \cite{EgorKond_1989}, \cite{RammSmir_2001}, \cite{ Wang_Hon_Ch_2006},
\cite{Qian&Fu&Xiong&Wei_2006}, \cite{Nakamura&Wang&Wang_2008} and \cite{SSS_CMAM}. In particular, the results of
\cite{SSS_CMAM} have opened perspective for further investigation of numerical methods for recovery of high-order
derivatives. Namely, the main criterion of the method's efficiency has been taken as its ability to achieve the optimal
order of accuracy by using a minimal amount of discrete information. Note that these aspects of numerical
differentiation remain still insufficiently studied. The present paper continues the research of \cite{SSS_CMAM},
\cite{Sem_Sol_2021}, \cite{Lyashko2022}   and proposes a numerical method for recovering the high-order mixed derivatives of smooth bivariate
functions. The method is not only stable to small perturbations of the input data, but is also optimal in terms of
accuracy and quantity of involved Fourier-Legendre coefficients, and also has a simple numerical implementation.
Note that the approach to numerical differentiation proposed below can be widely used in solving various practical problems.
For example, the expected results can be applied to the analysis of well-known
multidimensional Savitzky-Golay method \cite{Sh2016}, \cite{LYHB}.

The article is organized as follows. In Section \ref{DescrProb},  the problem statement for optimizing numerical differentiation methods in the sense of the minimal radius of Galerkin information is given. Sections \ref{TML2} and \ref{TMC} describe a modification of the truncation method and establish its accuracy estimates in quadratic and uniform metrics, respectively. Section \ref{SecMinimalRad} is devoted to finding order estimates for the minimal radius of Galerkin information while establishing the optimality (in the power scale) of the method under consideration. In Section \ref{CompEx}, numerical experiments will be provided to confirm the effectiveness of the proposed method.

For the further presentation of the subject matter of this paper, we need the following notation and concepts.
Let $\{\varphi_k(t)\}_{k=0}^\infty$  be the system of Legendre polynomials  orthonormal on $[-1,1]$ as
$$
\varphi_k(t)=\sqrt{k+1/2}(2^kk!)^{-1}\frac{d^k}{dt^k}[(t^2-1)^k], \quad k=0,1,2,\ldots .
$$

By $L_2=L_2(Q)$   we mean space of square-summable  on $Q=[-1,1]^2$ functions $f(t,\tau)$ with inner  product
$$
\langle f, g\rangle=\int_{-1}^{1}\int_{-1}^{1}f(t,\tau)g(t,\tau)d \tau d t
$$
and norm generated by it  
$$
\|f\|_{L_2}^2=\sum_{k,j=0}^{\infty}|\langle f, \varphi_{k,j} \rangle|^2 < \infty ,
$$
where $$ \langle f, \varphi_{k,j}\rangle=\int_{-1}^{1}\int_{-1}^{1}f(t,\tau)\varphi_k(t)\varphi_j(\tau)d\tau dt, \quad
k,j=0,1,2,\ldots,
$$
are Fourier-Legendre coefficients of $f$. Moreover, let $C=C(Q)$ be the space of continuous bivariate functions on $Q$ equipped with standard uniform norm  and  $\ell_p$, $1\leq p\leq\infty$, be the space of numerical sequences
$\overline{x}=\{x_{k,j}\}_{k,j\in\mathbb{N}_0}$, $\mathbb{N}_0=\{0\}\bigcup\mathbb{N}$, such that the corresponding
relation
\begin{equation}\label{noise}
	\|\overline{x}\|_{\ell_p}  := \left\{
	\begin{array}{cl}
		\bigg(\sum\limits_{k,j\in\mathbb{N}_0} |x_{k,j}|^p\bigg)^{\frac{1}{p}} < \infty ,
		\ & 1\leq p<\infty ,
		\\\\
		\sup\limits_{k,j\in\mathbb{N}_0}  |x_{k,j}| < \infty ,
		\ & p=\infty ,
	\end{array}
	\right.
\end{equation}
is fulfilled.

We introduce a space of functions
$$
L_{s,2}^\mu(Q)=\{f\in L_2(Q): \quad \|f\|_{s,\mu}^s:=\sum_{k,j=0}^{\infty} ({\underline{k\cdot j}})^{s\mu}|\langle f,
\varphi_{k,j}\rangle|^s<\infty\},
$$
where $\mu>0$,\ $1\le s<\infty$,\ $\underline{k}=\max\{1,k\}$,\ $k=0,1,2,\dots$.

Hereinafter we will use the same notations both for the space and for a unite ball from this space, which we call a class of functions and denote as
 $L_{s,2}^{\mu} = L_{s,2}^{\mu}(Q) = \{f\in
L_{s,2}^{\mu}\!: \|f\|_{s,\mu} \leq 1\}$. What exactly is meant by $
L_{s,2}^{\mu}$, space or class, will be clear depending on the context in each case. It should be noted that
$L^{\mu}_{s,2}$ is a generalization of the class of bivariate functions with dominating mixed partial derivatives.

We represent a function $f(t,\tau)$ from $L_{s,2}^{\mu}$, $\mu\geq 2r+1/2-1/s$, as
$$
f(t,\tau) = \sum_{k,j=0}^{\infty}\langle f, \varphi_{k,j}\rangle \varphi_k(t)\varphi_j(\tau),
$$
and by its mixed derivative $f^{(r,r)}$, $r=1,2,\ldots$, we mean the following series
\begin{equation}\label{r_deriv}
	f^{(r,r)}(t,\tau) =  \sum_{k,j=r}^{\infty}\langle f, \varphi_{k,j}\rangle \varphi^{(r)}_k(t)\varphi^{(r)}_j(\tau).
\end{equation}
Assume that instead of the exact values of  the Fourier-Legendre coefficients $\langle f, \varphi_{k,j} \rangle$ only
some their perturbations are known with the error level $\delta$ in the metrics of $\ell_p$, $1\leq p\leq\infty$.
More accurately, we assume that there is a sequence of numbers $\overline{f^\delta}= \{\langle
f^\delta, \varphi_{k,j} \rangle\}_{k,j\in\mathbb{N}_0}$ such that for $\overline{\xi}=
\{\xi_{k,j}\}_{k,j\in\mathbb{N}_0}$, where $\xi_{k,j}=\langle f-f^\delta,\varphi_{k,j}\rangle$, and for some $1\leq
p\leq \infty$ the relation
\begin{equation}\label{perturbation}
	\|\overline{\xi}\|_{\ell_p} \leq \delta , \quad 0<\delta <1 ,
\end{equation}
is true.

The research of this work is devoted to the optimization of methods for recovering the derivative (\ref{r_deriv}) of
functions from  $L^{\mu}_{s,2}$. Further, we give a strict statement of the problem to be studied. In the
coordinate plane $[r,\infty)\times[r,\infty)$ we take an arbitrary bounded domain $\Omega$. By $\card(\Omega)$ we mean
the number of points that make up $\Omega$ and by the information vector $G(\Omega, \overline{f}^{\delta})\in
\mathbb{R}^{N}$, $\card (\Omega) = N$, we take the set of perturbed values of Fourier-Legendre coefficients
$\left\{\langle f^{\delta}, \varphi_{k,j} \rangle\right\}_{(k,j)\in\Omega}$.

Let  $X=L_{2}(Q)$ or $X=C(Q)$.
By numerical differentiation algorithm we mean any mapping $\psi^{(r,r)} = \psi^{(r,r)}(\Omega)$ that corresponds to the information vector
$G(\Omega, \overline{f}^{\delta})$ an element $\psi^{(r,r)}(G(\Omega, \overline{f}^{\delta})) \in X$,
which is taken as an approximation to the derivative (\ref{r_deriv}) of function $f$ from the class $L_{s,2}^{\mu}$.
We denote by  $\Psi(\Omega)$ the set of all algorithms $\psi^{(r,r)}(\Omega):\,\mathbb{R}^{N}\rightarrow X$, that use
the same information vector $G(\Omega, \overline{f}^{\delta})$.

We do not require, generally speaking,  either linearity or even stability for algorithms from $\Psi(\Omega)$.
The only condition for these algorithms
is to use an input information in the form of perturbed values of the
Fourier-Legendre coefficients
with indices from the domain $\Omega$  of the coordinate plane. Such a general understanding of the algorithm
is explained by the desire to consider the widest range of possible methods of numerical differentiation.

The error of the algorithm $\psi^{(r,r)}$ on the class $L^{\mu}_{s,2}$ is determined by the quantity
$$
\varepsilon_{\delta}(L^{\mu}_{s,2}, \psi^{(r,r)}(\Omega), X, \ell_p)
= \sup_{\substack{f\in L^{\mu}_{s,2}, \\ \|f\|_{s,\mu}\leq 1}}
\ \, \sup\limits_{\overline{f}^{\delta}: \, (\ref{perturbation})}
\| f^{(r,r)} - \psi^{(r,r)}(G(\Omega, \overline{f}^{\delta})) \|_X .
$$

The minimal radius of the Galerkin information for the problem of numerical differentiation on the class
$L^{\mu}_{s,2}$ is given by
$$
R^{(r,r)}_{N,\delta} (L^{\mu}_{s,2}, X, \ell_p) = \inf\limits_{\Omega: \, \card(\Omega)\leq N}
\ \, \inf\limits_{\psi^{(r,r)}\in\Psi(\Omega)} \varepsilon_{\delta}(L^{\mu}_{s,2}, \psi^{(r,r)}(\Omega), X, \ell_p) .
$$

The quantity $R^{(r,r)}_{N,\delta} (L^{\mu}_{s,2}, X, \ell_p)$ describes the minimal possible accuracy in the metric of
space $X$, which can be achieved by numerical differentiation of arbitrary function $f\in L_{s,2}^{\mu}$ , while using
not more than $N$ values of its Fourier-Legendre coefficients that are $\delta$-perturbed in the $\ell_p$ metric. Note
that the minimal radius of Galerkin information in the problem of recovering the first partial derivatives was studied
in \cite{Sol_Stas_UMZ2022}, the mixed derivatives  $f^{(2,2)}$ in \cite{Lyashko2022}   and for other types of ill-posed problems, similar studies were previously carried out in
\cite{PS1996, Mileiko_Solodkii_2014}. It should be added that the minimal radius characterizes the information
complexity of the considered problem and is traditionally studied within the framework of the IBC-theory (Information-Based Complexity Theory), the foundations of which are laid in monographs \cite{TrW} and \cite{TrWW}.

The goal of our research is to find order-optimal estimates (in the power scale) for  $R^{(r,r)}_{N,\delta}
(L^{\mu}_{s,2}, C, \ell_p)$ and $R^{(r,r)}_{N,\delta} (L^{\mu}_{s,2}, L_{2}, \ell_p)$.

Finally, we introduce the symbolic notation for inequality and equality in order. 
	For two positive quantities $a$ and  $b$, we write $a\preceq b$ if there exists a constant $c > 0$,  such that $a\leq c b$. 
	We will write $a\asymp b$ if  $a\preceq b$ and $b\preceq a$.

\section{Truncation method. Error estimate in the metric of  $L_2$ }\label{TML2}

Many approaches to numerical differentiation have been developed so far (see e.g. \cite{Cul71}, \cite{And84}, \cite{Qu96}, \cite{RammSmir_2001} and also see \cite{SSS_Rew2021} and the references therein). However, most known methods have their drawbacks. So, in particular, some of them guarantee satisfactory accuracy only in the case of exactly given input data. Other methods are quite complex for numerical implementation. Here we mean the representation of the problem to be solved in the form of an integral equation (for example, in the framework of the Tikhonov method and its various variations), as well as difficult-to-implement rules for determining the regularization parameters (see \cite{RammSmir_2001}). Thus, the problem of constructing efficient methods for numerical differentiation remains open.
Recently in \cite{SSS_CMAM}, \cite{Sol_Stas_UMZ2022} a concise numerical method called the truncation method has been proposed as a simple and optimal (in the sense of minimal radius of Galerkin information) approach to numerically differentiating bivariate functions. The essence of this method is to replace the Fourier series (\ref{r_deriv}) with a finite Fourier sum over perturbed data $\langle f^\delta, \varphi_{k,j} \rangle$. In the truncation method, to ensure the stability of the approximation and achieve the required order accuracy, it is necessary to correctly choose a discretization parameter, which here serves as the regularization parameter. More strictly, the regularization process in the method under consideration consists of matching the discretization parameter with the perturbation level $\delta$ of the input data, the smoothness parameter $\mu$, the order of the derivative $r$, and the values of $p,s$. The simplicity of implementation is the main advantage of this method. Moreover, in Section \ref{SecMinimalRad} it will be established that the proposed version (\ref{ModVer}) of the truncation method is optimal in terms of the minimal radius of  Galerkin information. In other words, this method achieves the best order of accuracy while using the smallest possible amount of perturbed input data  $\langle f^\delta, \varphi_{k,j} \rangle$.

In the case of an arbitrary bounded domain $\Omega$ of the coordinate plane $[r,\infty)\times[r,\infty)$, the
truncation method for differentiating bivariate functions  has the form
$$
\mathcal{D}_\Omega^{(r,r)} f^\delta(t,\tau)=\sum_{(k,j)\in \Omega} \langle f^\delta, \varphi_{k,j} \rangle
\varphi_k^{(r)}(t)\varphi_j^{(r)}(\tau).
$$
To increase the efficiency of the approach under study, we take a hyperbolic cross as the domain $\Omega$ of
the following form
$$
\Omega = \Gamma_n := \{ (k,j):\ k\cdot j \leq rn-1,\quad k,j=r,\ldots, n-1\}, \qquad \card(\Gamma_n)\asymp n \ln n.
$$
Then  the proposed version of the truncation method can be written 
\begin{equation} \label{ModVer}
	\mathcal{D}_n^{(r,r)} f^\delta(t,\tau) = \sum_{k,j\ge r,\ kj\leq rn-1 } \langle f^\delta, \varphi_{k,j}\rangle
	\varphi^{(r)}_k(t)\varphi^{(r)}_j(\tau).
\end{equation}

We note that the idea of a hyperbolic cross for the problem of numerical differentiation was used earlier in the papers
\cite{Sem_Sol_2021}, \cite{Sol_Stas_UMZ2022},  \cite{SSS_CMAM}, \cite{Lyashko2022} (for more details about the usage of a hyperbolic cross in
solving the other  ill-posed problems see \cite{Pereverzev_Computing_1995}, \cite{ErbSem2015}, \cite{Mileiko_Solodkii_2016_ApAn}, \cite{Mileiko_Solodkii_2017_UMJ}). 

\begin{remark}
	In \cite{SSS_CMAM}, the problem of recovering the derivatives  $f^{(r,r)}$ of periodic functions was considered, 
	when the perturbed values
	of the Fourier coefficients
	w.r.t. the trigonometric system are taken as input information.
	Unfortunately, it is impossible to automatically transfer the results from the periodic case to the non-periodic one.
	In particular, this is because the best approximation accuracy for derivatives of non-periodic functions has a worse order than the best accuracy for approximation
	of derivatives of periodic functions (cf. \cite{SSS_CMAM}).
	Therefore, in the non-periodic case, to construct optimal methods for numerical differentiation,
	a modification of the previous methodology and the development of new techniques are required.
	At the moment, we have already constructed optimal methods for recovering the derivatives
	$f^{(1,1)}$   \cite{Sem_Sol_2021}, $f^{(2,2)}$ \cite{Lyashko2022},  $f^{(1,0)}$ \cite{Sol_Stas_UMZ2022} in the non-periodic case.
	The results of \cite{Sem_Sol_2021}, \cite{Sol_Stas_UMZ2022}, and \cite{Lyashko2022} together with the results of this work,
	create the ground and prospects for the development of optimal methods for recovering derivatives
	of any order for non-periodic functions of any number of variables. \end{remark}

\begin{remark}

If in (\ref{ModVer}) we put $r=0$, 
 then the problem of numerical differentiation is transformed into the problem of numerical summation. Earlier, to solve the problem of numerical summation of univariate functions, the truncation method was proposed in \cite{MatPer2002}. Further, this approach was extended to the case of bivariate functions in \cite{SolSha2015}, \cite{Sol_Stas_JC2020}. It should be noted that in the last two works, the idea of a  hyperbolic cross was successfully implemented, which made it possible to significantly reduce the computational resources without loss of accuracy. In the framework of the present investigation, we are going to achieve the same effect for the problem of numerical differentiation. 

		\end{remark}

Let us investigate the approximation properties of the method (\ref{ModVer}). To this end, we write the error of the method (\ref{ModVer}) as
\begin{equation}\label{fullError}
	f^{(r,r)}(t,\tau)-\mathcal{D}_n^{(r,r)} f^\delta(t,\tau)= \left(f^{(r,r)}(t,\tau)-\mathcal{D}_n^{(r,r)}
	f(t,\tau)\right)+\left(\mathcal{D}_n^{(r,r)} f(t,\tau)-\mathcal{D}_n^{(r,r)} f^\delta(t,\tau)\right).
\end{equation}
For the first difference on the right-hand side of (\ref{fullError}), the representation
\begin{equation}\label{Bound_err}
	f^{(r,r)}(t,\tau)-\mathcal{D}_n^{(r,r)} f(t,\tau)= \triangle_{1}(t,\tau)+\triangle_{2}(t,\tau)+\triangle_{3}(t,\tau)
\end{equation}
holds, where
\begin{equation}\label{Triangle_1HC}
	\triangle_{1}(t,\tau)= \sum_{k=n+1}^{\infty} \sum_{j=r}^{\infty} \langle f,
	\varphi_{k,j}\rangle\varphi^{(r)}_k(t)\varphi^{(r)}_j(\tau),
\end{equation}
\begin{equation}\label{Triangle_2HC}
	\triangle_{2}(t,\tau)= \sum_{k=r}^{n} \sum_{j=n+1}^{\infty} \langle f, \varphi_{k,j}\rangle\varphi^{(r)}_k(t)\varphi^{(r)}_j(\tau),
\end{equation}

\begin{equation}\label{Triangle_3HC}
	\triangle_{3}(t,\tau)= \sum_{k=r}^{n} \sum_{j=\frac{rn}{k}}^{n} \langle f,
	\varphi_{k,j}\rangle\varphi^{(r)}_k(t)\varphi^{(r)}_j(\tau).
\end{equation}

For our calculations, we need the following formula (see Lemma 18 \cite{Mul69})
\begin{equation}\label{Muller}
	\varphi_{k}'(t) \, = 2 \, \sqrt{k+1/2}
	\mathop{{\sum}^*}\limits_{l=0}^{k-1} \sqrt{l+1/2} \, \varphi_{l}(t) ,
	\ \ \ k\in\mathbb{N} ,
\end{equation}
where in aggregate \quad $ \mathop{{\sum}^*}\limits_{l=0}^{k-1} \sqrt{l+1/2} \, \varphi_{l}(t) \, $ \quad the summation
is extended over only those terms for which $k+l$ is odd.

In the sequel, we adopt the convention that $c$ denotes a generic positive coefficient, which
can vary from inequality to inequality and may only depend on basic parameters such as $\mu, r, p, s$
and others which may appear below. 

Let us estimate the error of the method (\ref{ModVer}) in the metric of $L_2$. An upper bound for difference
(\ref{Bound_err}) is contained in the following statement.

\begin{lemma}\label{lemma_BoundErrHC}
	Let $f\in L^\mu_{s,2}$, $1\leq s< \infty$, $\mu>2r+1/2-1/s$. Then it holds
	$$
	\|f^{(r,r)}-\mathcal{D}_n^{(r,r)} f\|_{L_2}\leq c\|f\|_{s,\mu} n^{-\mu+2r+1/2-1/s} \ln^{3/2-1/s} n.
	$$
\end{lemma}

\textit{Proof.} Using the formula (\ref{Muller}), from (\ref{Triangle_1HC}) we have
$$
\triangle_{1}(t,\tau)
= 4^r \sum_{k=n+1}^{\infty} \sum_{j=r}^{\infty} \sqrt{k+1/2}\sqrt{j+1/2}\:\langle
f, \varphi_{k,j}\:\rangle
$$
$$
\times \mathop{{\sum}^*}\limits_{l_1=r-1}^{k-1}  (l_1+1/2)\: \mathop{{\sum}^*}\limits_{l_2=r-2}^{l_1-1} (l_2+1/2)\ldots
\mathop{{\sum}^*}\limits_{l_{r-1}=1}^{l_{r-2}-1} (l_{r-1}+1/2)\mathop{{\sum}^*}\limits_{l_{r}=0}^{l_{r-1}-1}
\sqrt{l_{r}+1/2} \,\varphi_{l_r}(t) $$
$$
\times \mathop{{\sum}^*}\limits_{m_1=r-1}^{j-1} (m_1+1/2) \mathop{{\sum}^*}\limits_{m_2=r-2}^{m_1-1} (m_2+1/2)\ldots
\mathop{{\sum}^*}\limits_{m_{r-1}=1}^{m_{r-2}-1}  (m_{r-1}+1/2) \mathop{{\sum}^*}\limits_{m_{r}=0}^{m_{r-1}-1}  \sqrt{m_r+1/2}\:
\varphi_{m_r}(\tau) .
$$

We note that  in the representation $\triangle_{1}$ only those terms  take place for which all  indexes
$l_1+k, l_2+l_1,...,l_r+l_{r-1}, m_1+j, m_2+m_1,...,m_r+m_{r-1}$ are odd. Such rule is valid also for other terms, namely $\triangle_{2}, \triangle_{3}$ and
$\mathcal{D}_n^{(r,r)} f-\mathcal{D}_n^{(r,r)} f^\delta$, appearing in the error representation  (see (\ref{fullError}) -- (\ref{Triangle_3HC})).
In the following, for simplicity, we will omit the symbol "*"\, when denoting such summation operations, while taking into account this rule in the calculations.

Further, we change the order of summation and get
$$
\triangle_{1}(t,\tau)= \triangle_{11}(t,\tau)+\triangle_{12}(t,\tau),
$$
where
$$
\triangle_{11}(t,\tau)=4^r \mathop{{\sum}}\limits_{l_r=0}^{n-r+1}  \sqrt{l_r+1/2} \: \varphi_{l_r}(t)
\mathop{{\sum}}\limits_{m_r=0}^{\infty}  \sqrt{m_r+1/2} \: \varphi_{m_r}(\tau)
$$
\begin{equation}\label{triangle_{11}}
	\times \sum_{k=n+1}^{\infty} \sum_{j=m_r+r}^{\infty} \sqrt{k+1/2}\sqrt{j+1/2}\:\langle f, \varphi_{k,j}\:
	\rangle B^r_{k,j} ,
\end{equation}
$$
\triangle_{12}(t,\tau)=4^r \mathop{{\sum}}\limits_{l_r=n-r+2}^{\infty}  \sqrt{l_r+1/2}\: \varphi_{l_r}(t)
\mathop{{\sum}}\limits_{m_r=0}^{\infty}  \sqrt{m_r+1/2} \: \varphi_{m_r}(\tau)
$$
\begin{equation}\label{triangle_{12}}
	\times \sum_{k=l_r+r}^{\infty}\ \sum_{j=m_r+r}^{\infty} \sqrt{k+1/2}\sqrt{j+1/2}\:\langle f, \varphi_{k,j}\:
	\rangle B^r_{k,j}
\end{equation}
and
$$
B^r_{k,j}:=\mathop{{\sum}}\limits_{l_1=l_r+r-1}^{k-1}  (l_1+1/2)\mathop{{\sum}}\limits_{l_2=l_r+r-2}^{l_1-1}
(l_2+1/2) \ldots \mathop{{\sum}}\limits_{l_{r-1}=l_r+1}^{l_{r-2}-1}  (l_{r-1}+1/2)
$$
$$
\times \mathop{{\sum}}\limits_{m_1=m_r+r-1}^{j-1} (m_1+1/2) \mathop{{\sum}}\limits_{m_2=m_r+r-2}^{m_1-1} (m_2+1/2)\ldots
\mathop{{\sum}}\limits_{m_{r-1}=m_r+1}^{m_{r-2}-1} (m_{r-1}+1/2)
$$
\begin{equation}\label{axular1}
	\leq c (kj)^{2(r-1)} .
\end{equation}

At first, we consider the case $1<s<\infty.$
For $\triangle_{11}$ we have
$$
\|\triangle_{11}\|_{L_2}^2\leq 4^{2r} \mathop{{\sum}}\limits_{l_r=0}^{n-r+1}  (l_r+1/2)
\mathop{{\sum}}\limits_{m_r=0}^{\infty}  (m_r+1/2) \left(\sum_{k=n+1}^{\infty}\ \sum_{j=m_r+r}^{\infty} k^{\mu}
j^{\mu}\:|\langle f, \varphi_{k,j}\: \rangle| \frac{B_{k,j}^r} {(kj)^{\mu-1/2}}\right)^2 .
$$
Using H\"{o}lder inequality  and (\ref{axular1}),  for $\mu>2r+\frac{s-1}{s}-1/2$ we get
$$
\|\triangle_{11}\|_{L_2}^2\leq c \mathop{{\sum}}\limits_{l_r=0}^{n-r+1}  (l_r+1/2)
\mathop{{\sum}}\limits_{m_r=0}^{\infty}  (m_r+1/2)  \left(\sum_{k=n+1}^{\infty}\, \sum_{j=m_r+r}^{\infty}
k^{s\mu}j^{s\mu}\:|\langle f, \varphi_{k,j}\: \rangle|^s\right)^{2/s}
$$
$$
\times \left(\sum_{k=n+1}^{\infty}\, \sum_{j=m_r+r}^{\infty} (kj)^{-(\mu-2r+3/2)s/(s-1)}\right)^{2(s-1)/s}
$$
$$
\leq c  \|f\|_{s,\mu}^2 n^{-2(\mu-2r+3/2)+\frac{2(s-1)}{s}}  \mathop{{\sum}}\limits_{l_r=0}^{n-r+1}  (l_r+1/2)
\mathop{{\sum}}\limits_{m_r=0}^{\infty}  (m_r+1/2)^{-2(\mu-2r+3/2)+\frac{2(s-1)}{s}+1}
$$
$$
\leq c \|f\|_{s,\mu}^2 n^{-2(\mu-2r+3/2)+\frac{2(s-1)}{s}+2}.
$$
Thus, we find
$$
\|\triangle_{11}\|_{L_2}\leq c \|f\|_{s,\mu} n^{-\mu+2r+\frac{s-1}{s}-1/2} .
$$
Applying the estimating technique above we can bound the norm of $\triangle_{12}(t,\tau):$
$$
\|\triangle_{12}\|_{L_2}^2\leq c \|f\|_{s,\mu}^2 \mathop{{\sum}}\limits_{l_r=n-r+2}^{\infty}  (l_r+1/2)
\mathop{{\sum}}\limits_{m_r=0}^{\infty}  (m_r+1/2) \left(\sum_{k=l_r+r}^{\infty} \sum_{j=m_r+r}^{\infty}
(kj)^{-(\mu-2r+3/2)s/(s-1)}\right)^{2(s-1)/s}
$$
$$
\leq c \|f\|_{s,\mu}^2\mathop{{\sum}}\limits_{l_r=n-r+2}^{\infty}  (l_r+1/2)^{-2(\mu-2r+3/2)+\frac{2(s-1)}{s}+1}
\mathop{{\sum}}\limits_{m_r=0}^{\infty} (m_r+1/2)^{-2(\mu-2r+3/2)+\frac{2(s-1)}{s}+1}
$$
$$
\leq c \|f\|_{s,\mu}^2n^{-2(\mu-2r+3/2)+\frac{2(s-1)}{s}+2}.
$$
Summing up the estimates for $\triangle_{11}(t,\tau)$ and $\triangle_{12}(t,\tau)$ we obtain
$$
\|\triangle_{1}\|_{L_2}\leq\|\triangle_{11}\|_{L_2}+\|\triangle_{12}\|_{L_2}\leq c \|f\|_{s,\mu}
n^{-\mu+2r+\frac{s-1}{s}-1/2}.
$$

Now using the formula (\ref{Muller}), from (\ref{Triangle_2HC}) we have
$$
\triangle_{2}(t,\tau)= 4^r \sum_{k=r}^{n}\ \sum_{j=n+1}^{\infty} \sqrt{k+1/2}\sqrt{j+1/2}\:\langle f, \varphi_{k,j}
\rangle
$$
$$
\times \mathop{{\sum}}\limits_{l_1=r-1}^{k-1}  (l_1+1/2)\mathop{{\sum}}\limits_{l_2=r-2}^{l_1-1}  (l_2+1/2) \ldots
\mathop{{\sum}}\limits_{l_{r-1}=1}^{l_{r-2}-1}  (l_{r-1}+1/2) \mathop{{\sum}}\limits_{l_{r}=0}^{l_{r-1}-1}
\sqrt{l_r+1/2}\: \varphi_{l_r}(t)
$$
$$
\times \mathop{{\sum}}\limits_{m_1=r-1}^{j-1}  (m_1+1/2)\mathop{{\sum}}\limits_{m_2=r-2}^{m_1-1}  (m_2+1/2) \ldots
\mathop{{\sum}}\limits_{m_{r-1}=1}^{m_{r-2}-1}  (m_{r-1}+1/2)\mathop{{\sum}}\limits_{m_{r}=0}^{m_{r-1}-1}  \sqrt{m_r+1/2}\:
\varphi_{m_r}(\tau) .
$$
Further, we change the order of summation and get
$$
\triangle_{2}(t,\tau)= \triangle_{21}(t,\tau)+\triangle_{22}(t,\tau),
$$
where
$$
\triangle_{21}(t,\tau)=4^r \mathop{{\sum}}\limits_{l_r=0}^{n-r}  \sqrt{l_r+1/2} \, \varphi_{l_r}(t)
\mathop{{\sum}}\limits_{m_r=0}^{n-r+1}  \sqrt{m_r+1/2} \: \varphi_{m_r}(\tau)
$$
$$
\times \sum_{k=l_r+r}^{n}\ \sum_{j=n+1}^{\infty} \sqrt{k+1/2}\sqrt{j+1/2}\:\langle f, \varphi_{k,j}\rangle\: B^r_{k,j},
$$
$$
\triangle_{22}(t,\tau)=4^r \mathop{{\sum}}\limits_{l_r=0}^{n-r}  \sqrt{l_r+1/2} \, \varphi_{l_r}(t)
\mathop{{\sum}}\limits_{m_r=n-r+2}^{\infty}  \sqrt{m_r+1/2} \: \varphi_{m_r}(\tau)
$$
$$
\times \sum_{k=l_r+r}^{n} \sum_{j=m_r+r}^{\infty} \sqrt{k+1/2}\sqrt{j+1/2}\:\langle f, \varphi_{k,j}\:
\rangle B^r_{k,j}.
$$
Further, we estimate
$$
\|\triangle_{21}\|_{L_2}^2\leq 4^{2r} \mathop{{\sum}}\limits_{l_r=0}^{n-r}  (l_r+1/2)
\mathop{{\sum}}\limits_{m_r=0}^{n-r+1}  (m_r+1/2) \left( \sum_{k=l_r+r}^{n}\ \sum_{j=n+1}^{\infty} k^\mu j^\mu |\langle
f, \varphi_{k,j}\rangle|\frac{B^r_{k,j}}{(kj)^{\mu-1/2}}\right)^2 .
$$
Again, applying the H\"{o}lder inequality, we get for $\mu>2r+\frac{s-1}{s}-1/2$:
$$
\|\triangle_{21}\|_{L_2}^2
\leq c\|f\|_{s,\mu}^2n^{-2(\mu-2r+3/2)+\frac{2(s-1)}{s}}
$$
$$
\times \mathop{{\sum}}\limits_{l_r=0}^{n-r}  (l_r+1/2)^{-2(\mu-2r+3/2)+\frac{2(s-1)}{s}+1}
\mathop{{\sum}}\limits_{m_r=0}^{n-r+1}  (m_r+1/2)\leq c  \|f\|_{s,\mu}^2n^{-2(\mu-2r+3/2)+\frac{2(s-1)}{s}+2}.
$$
Further, we bound the norm of $\triangle_{22}:$
$$
\|\triangle_{22}\|_{L_2}^2\leq c \|f\|_{s,\mu}^2 \mathop{{\sum}}\limits_{l_r=0}^{n-r}  (l_r+1/2)
\mathop{{\sum}}\limits_{m_r=n-r+2}^{\infty}  (m_r+1/2)
\left(\sum_{k=l_r+r}^{n} \sum_{j=m_r+r}^{\infty}  (kj)^{-(\mu-2r+3/2)s/(s-1)}\right)^{2(s-1)/s}
$$
$$
\leq c \|f\|_{s,\mu}^2 n^{-2(\mu-2r+3/2)+\frac{2(s-1)}{s}+2} .
$$
Summing up the estimates for $\triangle_{21}$ and $\triangle_{22}$ we obtain
$$
\|\triangle_{2}\|_{L_2}\leq\|\triangle_{21}\|_{L_2}+\|\triangle_{22}\|_{L_2}\leq c\,
\|f\|_{s,\mu}n^{-\mu+2r+\frac{s-1}{s}-1/2}.
$$

Using the formula (\ref{Muller}), from (\ref{Triangle_3HC}) we have
$$
\triangle_{3}(t,\tau)= 4^r \sum_{k=r}^{n} \sum_{j=\frac{rn}{k}}^{n} \sqrt{k+1/2}\sqrt{j+1/2}\:\langle f,
\varphi_{k,j}\:\rangle
$$
$$
\times \mathop{{\sum}}\limits_{l_1=r-1}^{k-1} (l_1+1/2)\mathop{{\sum}}\limits_{l_2=r-2}^{l_1-1} (l_2+1/2)\ldots
\mathop{{\sum}}\limits_{l_{r-1}=1}^{l_{r-2}-1} (l_{r-1}+1/2) \mathop{{\sum}}\limits_{l_{r}=0}^{l_{r-1}-1}  \sqrt{l_r+1/2}\: \varphi_{l_r}(t)
$$
$$
\times \mathop{{\sum}}\limits_{m_1=r-1}^{j-1} (m_1+1/2)\mathop{{\sum}}\limits_{m_2=r-2}^{m_1-1} (m_2+1/2)\ldots
\mathop{{\sum}}\limits_{m_{r-1}=1}^{m_{r-2}-1} (m_{r-1}+1/2) \mathop{{\sum}}\limits_{m_{r}=0}^{m_{r-1}-1}  \sqrt{m_r+1/2}\: \varphi_{m_r}(\tau).
$$
Further, using arguments similar to those above, we get
$$
\|\triangle_{3}\|_{L_2}\leq c\, \|f\|_{s,\mu} n^{-\mu+2r+\frac{s-1}{s}-1/2}\ln^{3/2-1/s}n.
$$
The combination of (\ref{Bound_err}) and bounds for the norms of $\triangle_{1}$, $\triangle_{2}$, $\triangle_{3}$
makes it possible to establish the desired inequality.

In the case of $s=1$, Lemma is proved similarly. \vspace{0.1in}

\vskip -4mm

${}$ \ \ \ \ \ \ \ \ \ \ \ \ \ \ \ \ \ \ \ \ \ \ \ \ \ \ \ \ \ \ \ \ \ \ \ \ \ \
\ \ \ \ \ \ \ \ \ \ \ \ \ \ \ \ \ \ \ \ \ \
\ \ \ \ \ \ \ \ \ \ \ \ \ \ \ \ \ \ \ \ \ \
\ \ \ \ \ \ \ \ \ \ \ \ \ \ \ \ \ \ \ \ \ \
\ \ \ \ \ \ \ \ \ \ \ \ \ \ \ \ \ \ \ \ \ \
\ \ \ \ \ \ \ \ \ \ \ \ \ \ \ \ \ \ \ \ \ \
\ \ $\Box$

The following statement contains an estimate for the second difference from the right-hand side of (\ref{fullError}) in
the metric of $L_2$.

\begin{lemma}\label{lemma_BoundPertHC}
	Let condition (\ref{perturbation}) be satisfied.  Then for an arbitrary function $f\in L_2(Q)$ it holds
	$$
	\|\mathcal{D}^{(r,r)}_n f - \mathcal{D}^{(r,r)}_n f^\delta\|_{L_2} \leq c\delta n^{2r-1/p+1/2} \ln^{3/2-1/p} n .
	$$
\end{lemma}
\textit{Proof.} Let us write down the representation
$$\mathcal{D}^{(r,r)}_n f(t,\tau) - \mathcal{D}^{(r,r)}_n f^\delta(t,\tau)=
\sum_{k,j\ge r,\ kj\leq rn-1 } \langle f-f^\delta, \varphi_{k,j}\rangle \varphi^{(r)}_k(t)\varphi^{(r)}_j(\tau). $$
Using the formula (\ref{Muller}), we get
$$\mathcal{D}^{(r,r)}_n f(t,\tau) - \mathcal{D}^{(r,r)}_n f^\delta(t,\tau)= 4^r \sum_{k=r}^{n-1}\,
\sum_{j=r}^{\frac{rn-1}{k}} \sqrt{k+1/2}\sqrt{j+1/2}\:\langle f-f^\delta, \varphi_{k,j}\: \rangle
$$
$$
\times \mathop{{\sum}}\limits_{l_1=r-1}^{k-1}  (l_1+1/2)\: \mathop{{\sum}}\limits_{l_2=r-2}^{l_1-1} (l_2+1/2)\:
\ldots \mathop{{\sum}}\limits_{l_{r-1}=1}^{l_{r-2}-1} (l_{r-1}+1/2)\:
\mathop{{\sum}}\limits_{l_{r}=0}^{l_{r-1}-1} \sqrt{l_{r}+1/2}\: \varphi_{l_r}(t)
$$
$$
\times \mathop{{\sum}}\limits_{m_1=r-1}^{j-1}  (m_1+1/2)\: \mathop{{\sum}}\limits_{m_2=r-2}^{m_1-1} (m_2+1/2)\:
\ldots \mathop{{\sum}}\limits_{m_{r-1}=1}^{m_{r-2}-1} (m_{r-1}+1/2)\:
\mathop{{\sum}}\limits_{m_{r}=0}^{m_{r-1}-1} \sqrt{m_{r}+1/2}\: \varphi_{m_r}(\tau) .
$$
Further, we change the order of summation and get
$$
\mathcal{D}^{(r,r)}_n f(t,\tau) - \mathcal{D}^{(r,r)}_n f^\delta(t,\tau) = 4^r\,
\mathop{{\sum}}\limits_{l_r=0}^{n-r-1} \sqrt{l_r+1/2}\, \varphi_{l_r}(t)
\mathop{{\sum}}\limits_{m_r=0}^{\frac{rn-1}{l_r+r}-r}  \sqrt{m_r+1/2} \: \varphi_{m_r}(\tau)
$$
$$
\times \sum_{k=l_r+r}^{\frac{rn-1}{m_r+r}} \sum_{j=m_r+r}^{\frac{rn-1}{k}} \sqrt{k+1/2}\sqrt{j+1/2}\, \langle
f-f^\delta, \varphi_{k,j}\: \rangle B^r_{k,j} .
$$
Let $1<p<\infty$ first. Then, using the H\"{o}lder inequality and  estimate (\ref{axular1}), we find
$$
\|\mathcal{D}^{(r,r)}_n f - \mathcal{D}^{(r,r)}_n f^\delta\|_{L_2}^2
\leq
c \delta^2 \mathop{{\sum}}\limits_{l_r=0}^{n-r-1} (l_r+1/2)
\mathop{{\sum}}\limits_{m_r=0}^{\frac{rn-1}{l_r+r}-r}  (m_r+1/2)
\left(n^{\frac{(4r-3)p}{2(p-1)}+1} \sum_{k=l_r+r}^{\frac{rn-1}{m_r+r}} \frac{1}{k}\right)^{2(p-1)/p}
$$
$$
\leq c\delta^2 n^{4r-3+\frac{2(p-1)}{p}} \ln^{\frac{2(p-1)}{p}} n \mathop{{\sum}}\limits_{l_r=0}^{n-r-1} (l_r+1/2)
\mathop{{\sum}}\limits_{m_r=0}^{\frac{rn-1}{l_r+r}-r} (m_r+1/2) \asymp \delta^2 n^{4r-1+2(p-1)/p} \ln^{3-2/p} n ,
$$
which was required to prove.

In the cases of $p=1$ and $p=\infty$, Lemma is proved similarly. \vspace{0.1in}

\vskip -4mm

${}$ \ \ \ \ \ \ \ \ \ \ \ \ \ \ \ \ \ \ \ \ \ \ \ \ \ \ \ \ \ \ \ \ \ \ \ \ \ \
\ \ \ \ \ \ \ \ \ \ \ \ \ \ \ \ \ \ \ \ \ \
\ \ \ \ \ \ \ \ \ \ \ \ \ \ \ \ \ \ \ \ \ \
\ \ \ \ \ \ \ \ \ \ \ \ \ \ \ \ \ \ \ \ \ \
\ \ \ \ \ \ \ \ \ \ \ \ \ \ \ \ \ \ \ \ \ \
\ \ \ \ \ \ \ \ \ \ \ \ \ \ \ \ \ \ \ \ \ \
\ \ $\Box$

The combination of Lemmas \ref{lemma_BoundErrHC} and \ref{lemma_BoundPertHC} gives

\begin{theorem} \label{Th1}
	Let $f\in L^\mu_{s,2}$, $1\leq s< \infty$, $\mu>2r-1/s+1/2$,  condition (\ref{perturbation}) be satisfied.
	Then for $n\asymp \left(\delta^{-1} \ln^{1/p-1/s} \frac{1}{\delta}\right)^{\frac{1}{\mu-1/p+1/s}}$  the following bound is valid
	$$
	\|f^{(r,r)} - \mathcal{D}^{(r,r)}_n f^\delta\|_{L_2} \leq c \left(\delta \ln^{1/s-1/p}
	\frac{1}{\delta}\right)^{\frac{\mu-2r+1/s-1/2}{\mu-1/p+1/s}} \ln^{3/2-1/s} \frac{1}{\delta} .
	$$
\end{theorem}

\textit{Proof.}  Taking into account  Lemmas \ref{lemma_BoundErrHC}, \ref{lemma_BoundPertHC},  from (\ref{fullError})  we get
$$
\|f^{(r,r)} - \mathcal{D}^{(r,r)}_n f^\delta\|_{L_2} \leq \|f^{(r,r)}-\mathcal{D}_n^{(r,r)} f\|_{L_2}+ 	\|\mathcal{D}^{(r,r)}_n f - \mathcal{D}^{(r,r)}_n f^\delta\|_{L_2} $$
$$\leq c\|f\|_{s,\mu} n^{-\mu+2r+1/2-1/s} \ln^{3/2-1/s} n+c\delta n^{2r-1/p+1/2} \ln^{3/2-1/p} n. $$
Substituting the rule $n\asymp \left(\delta^{-1} \ln^{1/p-1/s} \frac{1}{\delta}\right)^{\frac{1}{\mu-1/p+1/s}}$  into the relation above completely proves Theorem.
$\Box$
\vskip 2mm

\begin{corollary} \label{Cor1}
	\rm In the considered problem, the truncation method  $\mathcal{D}^{(r,r)}_{n}$ (\ref{ModVer})
	achieves (in the  $L_2$-metric) the accuracy 
	$$ O\Big(\left(\delta \ln^{1/s-1/p}
	\frac{1}{\delta}\right)^{\frac{\mu-2r+1/s-1/2}{\mu-1/p+1/s}} \ln^{3/2-1/s} \frac{1}{\delta}\Big)
	$$
	on the class $L^{\mu}_{s,2}$, $\mu>2r-1/s+1/2$, and requires
	$$
	\card(\Gamma_{n}) \asymp
	n\ \ln n \asymp \left(\delta^{-1} \ln^{\mu} \frac{1}{\delta}\right)^{\frac{1}{\mu-1/p+1/s}}
	$$
	perturbed Fourier-Legendre coefficients.
\end{corollary}

\vskip 2mm

\begin{remark} \label{Rem1}
	\rm Let us consider the standard version of the truncation method with $\Omega=\Box_n:= [r,n]\times [r,n]$.
	Similar to \cite{Sem_Sol_2021}  one can check that such an approach guarantees (in the  $L_2$-metric) the accuracy 
	$$
	O\left(\delta^{\frac{\mu-2r+1/s-1/2}{\mu+2r-2/p+1/s+1/2}}\right)$$  on the class $L^{\mu}_{s,2}$, $\mu>2r-1/s+1/2$ , and
	requires
	$$
	\card(\Box_{n}) \asymp
	n^2 \asymp \delta^{-{\frac{2}{\mu+2r-2/p+1/s+1/2}}}
	$$
	perturbed Fourier-Legendre coefficients.
The difference between these two versions of the truncation method (standard and  (\ref{ModVer})) is that within the method (\ref{ModVer}) we use only the Fourier-Legendre coefficients $\langle f^\delta, \varphi_{k,j} \rangle$  with indexes $(k,j)$ from  $\Gamma_n\subset \Box_n$.  
The proposed approach (\ref{ModVer}) allows not only to reduce the amount of input information from $O(n^2)$ (for the standard version) to $O(n\ln n)$ (for (\ref{ModVer})), but also to achieve higher accuracy.
	

\end{remark}

\vskip 2mm
\section{Truncation method. Error estimate in the metric of $C$}\label{TMC}

Now we have to bound the error of (\ref{ModVer}) in the metric of $C$. An upper estimate for the norm of the difference
(\ref{Bound_err}) is contained in the following statement.

\begin{lemma}\label{lemma_BoundErrHCC}
	Let $f\in L^\mu_{s,2}$, $1\leq s< \infty$, $\mu>2r-1/s+3/2$. Then 
	$$
	\|f^{(r,r)}-\mathcal{D}^{(r,r)}_n f\|_{C}\leq c\|f\|_{s,\mu} n^{-\mu+2r-1/s+3/2}\ln^{2-1/s} n.
	$$
\end{lemma}

\textit{Proof.} Let us start with the case $1< s< \infty$. Using (\ref{triangle_{11}}) and (\ref{axular1}), we get
$$
\|\triangle_{11}\|_{C}
\leq c \|f\|_{s,\mu}  \mathop{{\sum}}\limits_{l_r=0}^{n-r+1}  (l_r+1/2) \mathop{{\sum}}\limits_{m_r=0}^{\infty}
(m_r+1/2) \left(\sum_{k=n+1}^{\infty}\ \sum_{j=m_r+r}^{\infty} (kj)^{\frac{(-\mu+2r-3/2)s}{s-1}}\right)^{(s-1)/s}
$$

$$
\leq c \|f\|_{s,\mu}  n^{-\mu+2r+3/2-1/s} .
$$
Moreover, from (\ref{triangle_{12}}) it follows
$$
\|\triangle_{12}\|_{C}
\leq c \|f\|_{s,\mu}  \mathop{{\sum}}\limits_{l_r=n-r+2}^{\infty} (l_r+1/2) \mathop{{\sum}}\limits_{m_r=0}^{\infty}
(m_r+1/2) \left(\sum_{k=l_r+r}^{\infty}\ \sum_{j=m_r+r}^{\infty} (kj)^{\frac{(-\mu+2r-3/2)s}{s-1}}\right)^{(s-1)/s}
$$

$$
\leq c \|f\|_{s,\mu}  n^{-\mu+2r+3/2-1/s} .
$$
Thus, we get
$$
\|\triangle_{1}\|_{C} \leq c \|f\|_{s,\mu}  n^{-\mu+2r+3/2-1/s} .
$$
Similarly, we find
$$
\|\triangle_{2}\|_{C} \leq c \|f\|_{s,\mu}  n^{-\mu+2r+3/2-1/s} ,
$$

$$
\|\triangle_{3}\|_{C} \leq c \|f\|_{s,\mu}  n^{-\mu+2r+3/2-1/s} \ln^{2-1/s} n .
$$
Substituting estimates for the norms of $\triangle_{1}$, $\triangle_{2}$, $\triangle_{3}$ into the relation
(\ref{Bound_err}) allows to establish the desired inequality.

In the case of $s=1$, Lemma is proved similarly. \vspace{0.1in}

\vskip -4mm

${}$ \ \ \ \ \ \ \ \ \ \ \ \ \ \ \ \ \ \ \ \ \ \ \ \ \ \ \ \ \ \ \ \ \ \ \ \ \ \
\ \ \ \ \ \ \ \ \ \ \ \ \ \ \ \ \ \ \ \ \ \
\ \ \ \ \ \ \ \ \ \ \ \ \ \ \ \ \ \ \ \ \ \
\ \ \ \ \ \ \ \ \ \ \ \ \ \ \ \ \ \ \ \ \ \
\ \ \ \ \ \ \ \ \ \ \ \ \ \ \ \ \ \ \ \ \ \
\ \ \ \ \ \ \ \ \ \ \ \ \ \ \ \ \ \ \ \ \ \
\ \ $\Box$

\begin{remark} 
By the conditions of  Lemma \ref{lemma_BoundErrHCC} we have that  $\mu-2r >3/2-1/s\geq 1/2$.  Thus,  $f^{(r,r)}$ is continuous and the problem of estimating $f^{(r,r)}$ in the $C$-metric is well-defined.
\end{remark}

The following statement contains an estimate for the second difference from the right-hand side of (\ref{fullError}) in
the metric of $C$.

\begin{lemma}\label{lemma_BoundPertHCC}
Assume that condition (\ref{perturbation}) is satisfied.   Then for an arbitrary function $f\in C$ the following bound is valid
	$$
	\|\mathcal{D}^{(r,r)}_n f - \mathcal{D}^{(r,r)}_n f^\delta\|_{C} \leq c \delta n^{2r+3/2-1/p} \ln^{2-1/p} n .
	$$
\end{lemma}
\textit{Proof.} Let $1< p<\infty$ first. Then, using the H\"{o}lder inequality and the estimate (\ref{axular1}), we find
$$
\|\mathcal{D}^{(r,r)}_n f - \mathcal{D}^{(r,r)}_n f^\delta\|_{C}
\leq c\delta n^{2r-3/2+\frac{p-1}{p}} \ln^{\frac{p-1}{p}} n \mathop{{\sum}}\limits_{l_r=0}^{n-r-1} (l_r+1/2)
\mathop{{\sum}}\limits_{m_r=0}^{\frac{rn-1}{l_r+r}-r} (m_r+1/2)
$$
$$
\leq c\delta n^{2r+1/2+\frac{p-1}{p}} \ln^{\frac{p-1}{p}} n \mathop{{\sum}}\limits_{l_r=0}^{n-r-1} \frac{1}{l_r+1/2} \asymp  \delta n^{2r+3/2-1/p} \ln^{2-1/p} n ,
$$
which was required to prove.

In the cases of $p=1$ and $p=\infty$, Lemma is proved similarly. \vspace{0.1in}

\vskip -4mm

${}$ \ \ \ \ \ \ \ \ \ \ \ \ \ \ \ \ \ \ \ \ \ \ \ \ \ \ \ \ \ \ \ \ \ \ \ \ \ \
\ \ \ \ \ \ \ \ \ \ \ \ \ \ \ \ \ \ \ \ \ \
\ \ \ \ \ \ \ \ \ \ \ \ \ \ \ \ \ \ \ \ \ \
\ \ \ \ \ \ \ \ \ \ \ \ \ \ \ \ \ \ \ \ \ \
\ \ \ \ \ \ \ \ \ \ \ \ \ \ \ \ \ \ \ \ \ \
\ \ \ \ \ \ \ \ \ \ \ \ \ \ \ \ \ \ \ \ \ \
\ \ $\Box$

The combination of Lemmas \ref{lemma_BoundErrHCC} and \ref{lemma_BoundPertHCC} gives

\begin{theorem} \label{Th2}
	Let $f\in L^\mu_{s,2}$, $1\leq s< \infty$, $\mu>2r-1/s+3/2$, and  condition (\ref{perturbation})
	be satisfied. Then for $n\asymp \left(\delta^{-1} \ln^{1/p-1/s} \frac{1}{\delta}\right)^{\frac{1}{\mu-1/p+1/s}}$ the following bound is valid
	$$
	\|f^{(r,r)}-\mathcal{D}^{(r,r)}_n f^\delta\|_{C} \leq c \left(\delta \ln^{1/s-1/p}
	\frac{1}{\delta}\right)^{\frac{\mu-2r+1/s-3/2}{\mu-1/p+1/s}} \ln^{2-1/s} \frac{1}{\delta}  .
	$$
\end{theorem}

\vskip 2mm

\begin{corollary} \label{Cor2}
	\rm In the considered problem, the truncation method  $\mathcal{D}^{(r,r)}_{n}$ (\ref{ModVer})
	achieves (in the $C$-metric) the accuracy  
	$$
	O\Big(\left(\delta \ln^{1/s-1/p} \frac{1}{\delta}\right)^{\frac{\mu-2r+1/s-3/2}{\mu-1/p+1/s}} \ln^{2-1/s}
	\frac{1}{\delta}\Big)
	$$
	on the class $L^{\mu}_{s,2}$, $\mu>2r-1/s+3/2$, and requires
	$$
	\card(\Gamma_{n}) \asymp
	n \ln n \asymp\left(\delta^{-1} \ln^{\mu} \frac{1}{\delta}\right)^{\frac{1}{\mu-1/p+1/s}}
	$$
	perturbed Fourier-Legendre coefficients.
\end{corollary}

\vskip 2mm

\begin{remark} \label{Rem2}
	\rm Consider the standard variant of the truncation method with $\Omega=\Box_n$.
Similar to \cite{Sem_Sol_2021} one can verify that this approach guarantees (in the $C$-metric) the accuracy
	$$
	O\left(\delta^{\frac{\mu-2r+1/s-3/2}{\mu+2r-2/p+1/s+3/2}}\right)
	$$
	on the class $L^{\mu}_{s,2}$, $\mu>2r-1/s+3/2$, and requires
	$$
	\card(\Box_{n}) \asymp
	n^2 \asymp \delta^{-{\frac{2}{\mu+2r-2/p+1/s+3/2}}}
	$$
	perturbed Fourier-Legendre coefficients. Comparison of the estimates found above with the corresponding estimates for the method
	$\mathcal{D}^{(r,r)}_{n}$ (\ref{ModVer}) (see Corollary \ref{Cor2}) demonstrates that (\ref{ModVer}) is more efficient
	both in terms of accuracy and the amount of discrete information used.
\end{remark}

\vskip 2mm

\begin{remark} \label{Rem3} 
	\rm Earlier, the method $\mathcal{D}^{(r,r)}_{n}$ (\ref{ModVer}) was studied  for the problem of
	numerical differentiation of functions from
	$L^{\mu}_{s,2}$  in the case of $r=1$ \& $p=s=2$ (see \cite{Sem_Sol_2021}) and in the case of $r=2$ \& $s=2$ (see \cite{Lyashko2022}). Thus, the results of Theorems \ref{Th1} and \ref{Th2} generalize
	studies of \cite{Sem_Sol_2021}, \cite{Lyashko2022} for the case of arbitrary $r,p,s$. 
\end{remark}

\vskip 2mm
\section{Minimal radius of Galerkin information}\label{SecMinimalRad}
Now, we are in the position to find sharp estimates (in the power scale) for the minimal radius. First, we establish a lower
estimate for the quantity $R_{N,\delta}^{(r,r)}(L^{\mu}_{s,2}, C, \ell_p)$. We fix an arbitrarily chosen domain
$\hat{\Omega}$, $\card(\hat{\Omega})\leq N$, of the coordinate plane $[r,\infty)\times[r,\infty)$ and construct an auxiliary function
$$
f_1(t,\tau)  = \widetilde{c} \, \bigg( \varphi_0(t) \varphi_0(\tau) \,
+ \,  N^{-\mu-1/s}\, r^{-\mu} \varphi_r(\tau) \mathop{{\sum}\, '}\limits_{k=N+r}^{3N+r}
\varphi_k(t) \bigg) ,
$$
where the sum $\mathop{{\sum}\, '}\limits_{k=N+r}^{3N+r}$ is taken over any $N$ pairwise distinct functions $\varphi_k(t)$ such that
$N+r\leq k\leq 3N+r$ and $(k,r)\notin \hat{\Omega}$. There is at least one set of such functions.

Now we estimate the norm of $f_1$ in the space metric  $L^{\mu}_{s,2}$:
$$
\|f_1\|^s_{s,\mu} = {\widetilde{c}}^{\,s} \, \bigg( 1
+  N^{-s\mu-1}  \,  \mathop{{\sum}\, '}\limits_{k=N+r}^{3N+r}
k^{s\mu} \bigg) \leq
{\widetilde{c}}^{\,s}    \,  \bigg( 1
+ 4^{s\mu} \bigg) .
$$
Whence it follows that to satisfy the condition $\|f_1\|_{s,\mu}\leq 1$ it suffices to take
\begin{equation}\label{c_wave}
{\widetilde{c}} =  \bigg( 1
+ 4^{s\mu} \bigg)^{-1/s} .
\end{equation}

Next, we take another function from the class $L^{\mu}_{s,2}$:
$$
f_2(t,\tau) = \widetilde{c} \, \varphi_0(t) \varphi_0(\tau)  .
$$

Let us find a lower bound for the quantity $\|f_1^{(r,r)}-f_2^{(r,r)}\|_{C}$. For this we need formulas
$$
\varphi_r^{(r)}(t) = \frac{\sqrt{r+1/2}}{2^{r-1/2}}\ \frac{(2r)!}{r!}\ \varphi_0(t), \quad f_2^{(r,r)}(t,\tau) \equiv 0
,
$$
$$
f_1^{(r,r)}(t,\tau) = \frac{\widetilde{c}}{r^\mu} \, {N^{-\mu-1/s}} \,
\varphi_r^{(r)}(\tau) \mathop{{\sum}\, '}\limits_{k=N+r}^{3N+r} \varphi_k^{(r)}(t)
$$
$$
= 2^r\, \frac{\widetilde{c}}{r^\mu} \, {N^{-\mu-1/s}} \, \varphi_r^{(r)}(\tau) \mathop{{\sum}\,
	'}\limits_{k=N+r}^{3N+r} \sqrt{k+1/2} \mathop{{\sum}}\limits_{l_1=r-1}^{k-1} (l_1+1/2)
\mathop{{\sum}}\limits_{l_2=r-2}^{l_1-1} (l_2+1/2)
$$
\begin{equation}  \label{f_1^(r,r)}
	\ldots \mathop{{\sum}}\limits_{l_{r-1}=1}^{l_{r-2}-1} (l_{r-1}+1/2) \mathop{{\sum}}\limits_{l_{r}=0}^{l_{r-1}-1}
	\sqrt{l_r+1/2}\, \varphi_{l_r}(t) .
\end{equation}
We note that in the right-hand side of (\ref{f_1^(r,r)}) only terms with odd indexes $l_1+k, l_2+l_1,..., l_r+l_{r-1}$ take part.

It is easy to see that
$$
	\|f_1^{(r,r)}-f_2^{(r,r)}\|_{C}
	\geq  |f_1^{(r,r)}(1,1)|
	\geq  \overline{c} \, N^{-\mu+2r-1/s+3/2}  ,
$$
where
\begin{equation}\label{c_line}
\overline{c} = \frac{\widetilde{c}\sqrt{r+1/2}}{4^r r^\mu}\,  \frac{(2r)!}{(r!)^2} .
\end{equation}
Since for any $1\leq p\leq \infty$ it holds true
$$
\|\overline{f}_1-\overline{f}_2\|_{\ell_p}
= \frac{\widetilde{c}}{r^\mu} \, N^{-\mu-1/s+1/p} ,
$$
then in the case of $N^{-\mu-1/s+1/p}\leq r^\mu \delta/{\widetilde{c}}$  ,  the functions
$$
f^\delta_1(t,\tau) = f_2(t,\tau), \qquad f^\delta_2(t,\tau) = f_1(t,\tau)
$$
can be considered as   $\delta$-perturbations of $f_1 $ and $f_2$, respectively.

Let us find  an upper bound  for $\|f_1^{(r,r)}-f_2^{(r,r)}\|_{C}$.
Taking into account the relation $G(\hat{\Omega},\overline{f}_1^{\delta})=G(\hat{\Omega},\overline{f}_2^{\delta})$,
for any $\psi^{(r,r)}(\hat{\Omega})\in\Psi(\hat{\Omega})$ we find
$$
\|f_1^{(r,r)}-f_2^{(r,r)}\|_{C}
\leq \|f_1^{(r,r)}-\psi^{(r,r)}(G(\hat{\Omega},\overline{f}_1^{\delta}))\|_{C}
+ \|f_2^{(r,r)}-\psi^{(r,r)}(G(\hat{\Omega},\overline{f}_2^{\delta}))\|_{C}
\leq
$$
$$
\leq 2 \, \sup_{\substack{f\in L^{\mu}_{s,2},\\ \|f\|_{s,\mu}\leq 1}}
\ \sup_{\substack{\overline{f^\delta}: \, (\ref{perturbation})}}
\| f^{(r,r)} - \psi^{(r,r)}(G(\hat{\Omega},\overline{f}^{\delta})) \|_C
=: 2 \, \varepsilon_{\delta}(L^{\mu}_{s,2}, \psi^{(r,r)}(\hat{\Omega}), C, \ell_p) .
$$

That is
$$
\varepsilon_{\delta}(L^{\mu}_{s,2}, \psi^{(r,r)}(\hat{\Omega}), C, \ell_p)
\geq \frac{\overline{c}}{2} \, N^{-\mu+2r-1/s+3/2} .
$$

From the fact that the domain $\hat{\Omega}$
and the algorithm $\psi^{(r,r)}(\hat{\Omega})\in\Psi(\hat{\Omega})$  are arbitrary, follows that
$$
R_{N,\delta}^{(r,r)}(L^{\mu}_{s,2}, C, \ell_p)
\geq \frac{\overline{c}}{2} \, N^{-\mu+2r-1/s+3/2} .
$$

Thus, the following assertion is proved.

\begin{theorem} \label{Th5.1}
	Let $1\leq s< \infty$, $\mu>2r-1/s+3/2$, $1\leq p \leq \infty$, $N\geq\Big(r^\mu\delta/\widetilde{c}\Big)^{-1/(\mu+1/s-1/p)}$.  Then 
	$$
	R_{N,\delta}^{(r,r)}(L^{\mu}_{s,2}, C, \ell_p)
	\geq \frac{\overline{c}}{2} \, N^{-\mu+2r-1/s+3/2} ,
	$$
	where the constants $\widetilde{c}$ and $\overline{c}$ are defined by (\ref{c_wave}) and (\ref{c_line}) respectively.  
\end{theorem}

The following assertion contains sharp estimates (in the power scale) for the minimal radius in the uniform metric.

\begin{theorem} \label{Th5.2}
	Let $1\leq s< \infty$, $\mu>2r-1/s+3/2$, $1\leq p \leq \infty$.
	Then for $N \asymp \left(\delta^{-1} \ln^{\mu}
	\frac{1}{\delta}\right)^{\frac{1}{\mu-1/p+1/s}}$ it holds
	$$
	N^{-\mu+2r-1/s+3/2} \preceq
	R_{N,\delta}^{(r,r)}(L^{\mu}_{s,2}, C, \ell_p)
	\preceq N^{-\mu+2r-1/s+3/2} \ln^{\mu-2r+1/2} N
	$$
	or
	$$
	\left(\delta \ln^{-\mu} \frac{1}{\delta}\right)^{\frac{\mu-2r+1/s-3/2}{\mu-1/p+1/s}} \preceq
	R_{N,\delta}^{(r,r)}(L^{\mu}_{s,2}, C, \ell_p) \preceq \left(\delta \ln^{1/s-1/p}
	\frac{1}{\delta}\right)^{\frac{\mu-2r+1/s-3/2}{\mu-1/p+1/s}} \ln^{2-1/s} \frac{1}{\delta} .
	$$
	The upper bound is realized  by (\ref{ModVer}) for $n\asymp \left(\delta^{-1} \ln^{1/p-1/s}
	\frac{1}{\delta}\right)^{\frac{1}{\mu-1/p+1/s}}$ .
\end{theorem}

\bf Proof.
\rm
The upper bound for  $R_{N,\delta}^{(r,r)}(L^{\mu}_{s,2}, C, \ell_p)$ follows from Theorem \ref{Th2}.
The lower bound is found in Theorem  \ref{Th5.1}.

\vskip -4mm

${}$ \ \ \ \ \ \ \ \ \ \ \ \ \ \ \ \ \ \ \ \ \ \ \ \ \ \ \ \ \ \ \ \ \ \ \ \ \ \
\ \ \ \ \ \ \ \ \ \ \ \ \ \ \ \ \ \ \ \ \ \
\ \ \ \ \ \ \ \ \ \ \ \ \ \ \ \ \ \ \ \ \ \
\ \ \ \ \ \ \ \ \ \ \ \ \ \ \ \ \ \ \ \ \ \
\ \ \ \ \ \ \ \ \ \ \ \ \ \ \ \ \ \ \ \ \ \
\ \ \ \ \ \ \ \ \ \ \ \ \ \ \ \ \ \ \ \ \ \
\ \ $\Box$

\rm
Let us turn to estimate the minimal radius in the integral metric.

\begin{theorem} \label{Th5.4}
	Let $1\leq s< \infty$, $\mu>2r-1/s+1/2$, $1\leq p \leq \infty$. Then for any
	$N\geq\Big(r^\mu\delta/\widetilde{c}\Big)^{-1/(\mu+1/s-1/p)}$ it holds
	$$
	R_{N,\delta}^{(r,r)}(L^{\mu}_{s,2}, L_2, \ell_p) \geq \overline{\overline{c}} \, N^{-\mu+2r-1/s+1/2} ,
	$$
	where $\overline{\overline{c}} = \frac{\widetilde{c} \sqrt{r+1/2}}{2^{4r-3/2}\, r^{\mu-1}} \frac{(2r)!}{(r!)^2}$ and $\widetilde{c}$ is defined by (\ref{c_wave}).
\end{theorem}

\rm

\bf Proof. \rm The proof of Theorem \ref{Th5.4} almost completely coincides with the proof of Theorem \ref{Th5.1},
including the form of the auxiliary functions
$f_1$, $f_1^{\delta}$, $f_2$, $f_2^{\delta}$.
The only difference is in the lower estimate of the norm of the difference
$f_1^{(r,r)}-f_2^{(r,r)}$. Changing the order of summation in  (\ref{f_1^(r,r)}) yields to the representation
$$
f_1^{(r,r)}(t,\tau) = 2^r\, \frac{\widetilde{c}}{r^\mu} \, {N^{-\mu-1/s}} \, \varphi_r^{(r)}(\tau)
\Big(\mathop{{\sum}}\limits_{l_r=0}^{N} \sqrt{l_r+1/2}\, \varphi_{l_r}(t) \mathop{{\sum}\,
	'}\limits_{k=N+r}^{3N+r} \sqrt{k+1/2}
$$
$$
+ \mathop{{\sum}}\limits_{l_r=N+1}^{3N} \sqrt{l_r+1/2}\, \varphi_{l_r}(t) \mathop{{\sum}\, '}\limits_{k=l_r+r}^{3N+r}
\sqrt{k+1/2}\Big) \mathop{{\sum}}\limits_{l_1=l_r+r-1}^{k-1} (l_1+1/2) \mathop{{\sum}}\limits_{l_2=l_r+r-2}^{l_1-1}
(l_2+1/2) \ldots \mathop{{\sum}}\limits_{l_{r-1}=l_r+1}^{l_{r-2}-1} (l_{r-1}+1/2)  .
$$
It is easy to verify that
$$
\|f_1^{(r,r)}-f_2^{(r,r)}\|_{L_2}^2
\geq  \frac{c'}{4^{4r-2}((r-1)!)^2} N^{-2\mu+4r-2/s+1} ,
$$
where
$$
c' = \frac{2 \widetilde{c}^2}{r^{2\mu}} (r+1/2) \frac{((2r)!)^2}{(r!)^2} .
$$
Whence we obtain the relation
$$
\varepsilon_{\delta}(L^{\mu}_{s,2}, \psi^{(r,r)}(\hat{\Omega}), L_2, \ell_p)
\geq \overline{\overline{c}} \, N^{-\mu+2r-1/s+1/2}
$$
is true for any $N\geq\Big(r^\mu\delta/\widetilde{c}\Big)^{-1/(\mu+1/s-1/p)}$. From the fact that the domain
$\hat{\Omega}$ and the algorithm $\psi^{(r,r)}(\hat{\Omega})\in\Psi(\hat{\Omega})$  are arbitrary, it follows that
$$
R_{N,\delta}^{(r,r)}(L^{\mu}_{s,2}, L_2, \ell_p)
\geq \overline{\overline{c}} \, N^{-\mu+2r-1/s+1/2} .
$$
Thus, the proof of Theorem \ref{Th5.4} has been complied.

\vskip -4mm

${}$ \ \ \ \ \ \ \ \ \ \ \ \ \ \ \ \ \ \ \ \ \ \ \ \ \ \ \ \ \ \ \ \ \ \ \ \ \ \
\ \ \ \ \ \ \ \ \ \ \ \ \ \ \ \ \ \ \ \ \ \
\ \ \ \ \ \ \ \ \ \ \ \ \ \ \ \ \ \ \ \ \ \
\ \ \ \ \ \ \ \ \ \ \ \ \ \ \ \ \ \ \ \ \ \
\ \ \ \ \ \ \ \ \ \ \ \ \ \ \ \ \ \ \ \ \ \
\ \ \ \ \ \ \ \ \ \ \ \ \ \ \ \ \ \ \ \ \ \
\ \ $\Box$

The following statement contains sharp estimates (in the power scale) for the minimal radius in the integral metric.
\rm

\begin{theorem} \label{Th5.5}
	Let $1\leq s< \infty$, $\mu>2r-1/s+1/2$, $1\leq p \leq \infty$. Then for $N\asymp \Big(\delta^{-1}\ \ln^{\mu}
	\frac{1}{\delta}\Big)^{1/(\mu-1/p+1/s)}$ it holds
	$$
	N^{-\mu+2r-1/s+1/2} \preceq
	R_{N,\delta}^{(r,r)}(L^{\mu}_{s,2}, L_2, \ell_p)
	\preceq N^{-\mu+2r-1/s+1/2} \ln^{\mu-2r+1} N
	$$
	or
	$$
	\left(\delta \ln^{-\mu} \frac{1}{\delta}\right)^{\frac{\mu-2r+1/s-1/2}{\mu-1/p+1/s}} \preceq
	R_{N,\delta}^{(r,r)}(L^{\mu}_{s,2}, L_2, \ell_p) \preceq \left(\delta \ln^{1/s-1/p}
	\frac{1}{\delta}\right)^{\frac{\mu-2r+1/s-1/2}{\mu-1/p+1/s}} \ln^{3/2-1/s} \frac{1}{\delta} .
	$$
The upper bound is realized by (\ref{ModVer}) for  $n\asymp
	\left(\delta^{-1} \ln^{1/p-1/s} \frac{1}{\delta}\right)^{\frac{1}{\mu-1/p+1/s}}$.
\end{theorem}

{\textbf{\textit{Proof.}}
	\rm The upper bound for  $R_{N,\delta}^{(r,r)}(L^{\mu}_{s,2}, L_2, \ell_p)$ follows from Theorem
	\ref{Th1}.
	The lower bound is found in Theorem \ref{Th5.4}.
	
	\vskip -4mm
	
	${}$ \ \ \ \ \ \ \ \ \ \ \ \ \ \ \ \ \ \ \ \ \ \ \ \ \ \ \ \ \ \ \ \ \ \ \ \ \ \
	\ \ \ \ \ \ \ \ \ \ \ \ \ \ \ \ \ \ \ \ \ \
	\ \ \ \ \ \ \ \ \ \ \ \ \ \ \ \ \ \ \ \ \ \
	\ \ \ \ \ \ \ \ \ \ \ \ \ \ \ \ \ \ \ \ \ \
	\ \ \ \ \ \ \ \ \ \ \ \ \ \ \ \ \ \ \ \ \ \
	\ \ \ \ \ \ \ \ \ \ \ \ \ \ \ \ \ \ \ \ \ \
	\ \ $\Box$

	\section{Computational experiments}\label{CompEx}
	
	To demonstrate the effectiveness of the proposed method for recovering high-order derivatives some numerical experiments were carried out.
	The calculations were performed on a computer with a 4-core Intel Core i5 processor and 16 GB memory in the mathematical modeling environment MATLAB 2022a.
	

	\subsection{Example 1}
	
	We consider the function  $F_1(t,\tau) = f(t)f(\tau)/C$,  where $C=754$  and
	$$
	f(t)=\left\{
	\begin{array}{cl}
		-1/8 t^2+1/12 t^4-1/20 t^5+1/42 t^7-3/224 t^8,   & -1\leq t<0 ,
		\\\\
		-1/8 t^2+1/12 t^4-1/20 t^5+1/45 t^7-3/240 t^8, &  0\leq t \leq 1.
	\end{array}
	\right.
	$$	
	Taking into account the definition of space $L_{s,2}^\mu$ it is easy to see that $\|F_1\|_{2,\mu} \approx 1$ for $\mu=5,5$ and $\|F_1^{(2,2)}\|_{L_2}\approx 10^{-4}.$
	
	The simulation of the noise in the input data was done in two different ways:
	
	\begin{itemize}
		\item a random noise adds to the values of the Fourier-Legendre coefficients. The noise is generated by the
		$\mbox{{\bf randn}({\bf size}}({\cal F})) \delta$ command, where {\bf randn}
		and {\bf size} are standard functions of the MATLAB system, and ${\cal F}$ is a matrix for exact values of the Fourier-Legendre coefficients;
		
		\item the values of the Fourier-Legendre coefficients recovered by the quadrature trapezoid formula on a uniform grid with a step $h$   so that condition (\ref{perturbation}) is satisfied for a given $\delta$.

	\end{itemize}

	Numerical experiments were carried out for the following error levels: $\delta= 10^{-6}, 10^{-7}, 10^{-8}$. Tables \ref{tbl1} and \ref{tbl2} show the results of numerical calculations for the approximation of $F_1^{(2,2)}$ by the truncation method (2.1) for two different types of noise (random data noise and trapezoid formula errors).  "Error in $L_2$" and "Error in $C$" columns contain the recovery accuracy in the $L_2-$ and $C-$ metrics respectively, $n$ and  $\card(\Gamma_{n})$  columns indicate the highest degree of Legendre polynomial and the number of Fourier-Legendre coefficients involved, resp. Also in Table \ref{tbl2}, $h$ means the step size in the quadrature formula.
	
	The graphs \ref{Fig2} and \ref{Fig1} show the exact derivative $F_1{(2,2)}$ and its approximations constructed on data with random noise and with a noise, generated by the trapezoid formula, respectively.

	%

	\begin{table}[h!]
		\centering
		\caption{ The results of recovering derivative $F_1^{(2,2)}$ for random noise }
		\label{tbl1}
		\begin{tabular}{|c|c|c|c|}
			\hline
			$\delta$ & $10^{-6}$ & $10^{-7}$ & $10^{-8}$  \\ \hline
			Error in $L_2$       &  $1,1 \cdot 10^{-4} $       &  $2,73 \cdot 10^{-5} $  &  $6,7  \cdot 10^{-6} $  \\ \hline
			Error in $C$       &  $1,2 \cdot 10^{-3} $       &  $3,4 \cdot 10^{-4} $  &  $5 \cdot10^{-5} $  \\ \hline
			n  &  $19$      & $24$    & $31$    \\ \hline
			$\card(\Gamma_{n})$  &  $52$      &  $80$    & $106$    \\ \hline
		\end{tabular}
	\end{table}
	\begin{table}[ht]
		\centering
		\caption{ The results of recovering derivative $F_1^{(2,2)}$ for noise from quadrature formula }
		\label{tbl2}
		\begin{tabular}{|c|c|c|c|}
			\hline
			$\delta$ & $10^{-6}$ & $10^{-7}$ & $10^{-8}$  \\ \hline
			Error in $L_2$          &  $4,8 \cdot 10^{-5} $       &  $3,2 \cdot 10^{-5} $  &  $6,6  \cdot 10^{-6} $  \\ \hline
			Error in $C$        &  $7,53 \cdot 10^{-4} $       &  $4,9 \cdot 10^{-4} $  &  $2,53 \cdot10^{-5} $  \\ \hline
			n  &  $19$      & $24$    & $31$    \\ \hline
			h  &  $1,16\cdot 10^{-4}$    & $8\cdot10^{-5} $  & $ 4\cdot10^{-5}$    \\ \hline
		\end{tabular}
	\end{table}

	\begin{figure}[h!]
		\begin{minipage}[h]{0.5\linewidth}
			\center{\includegraphics[width=1\linewidth]{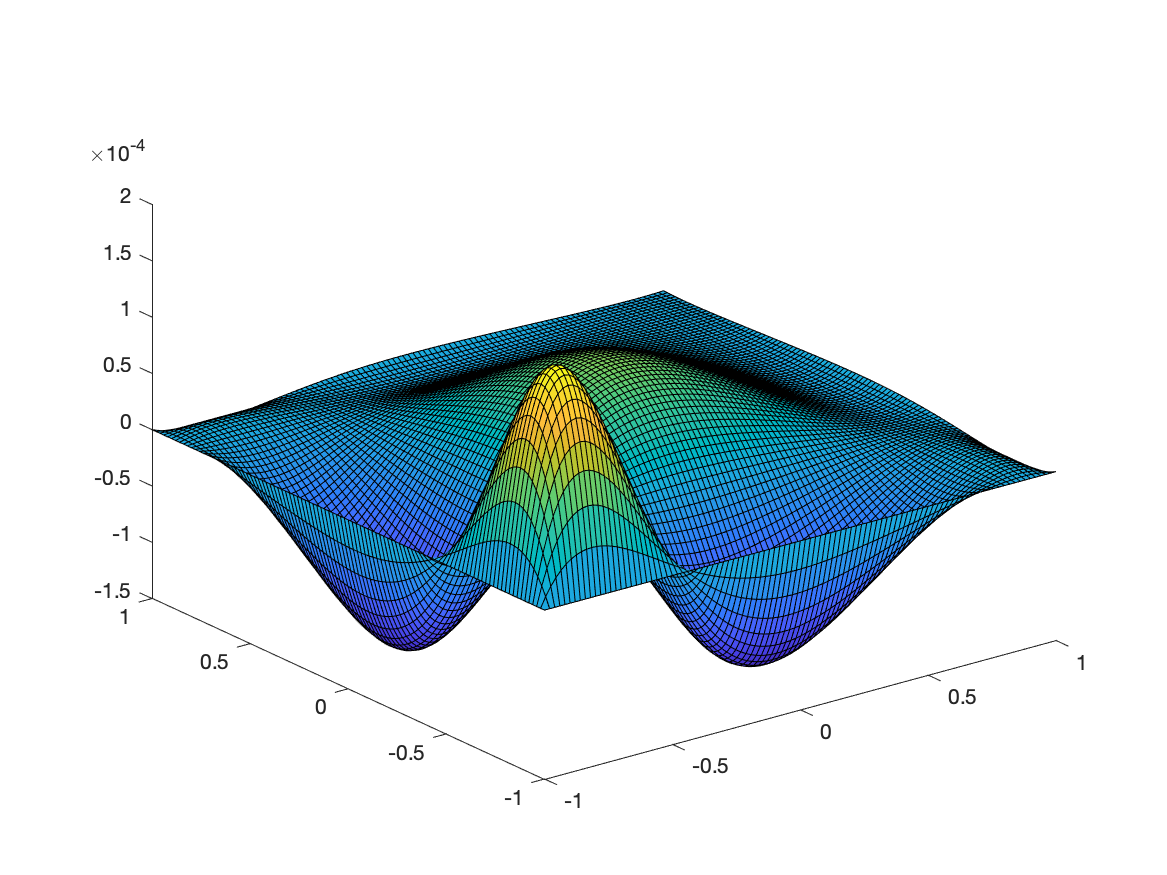} \\ a}
		\end{minipage}
		\begin{minipage}[h]{0.5\linewidth}
			\center{\includegraphics[width=1\linewidth]{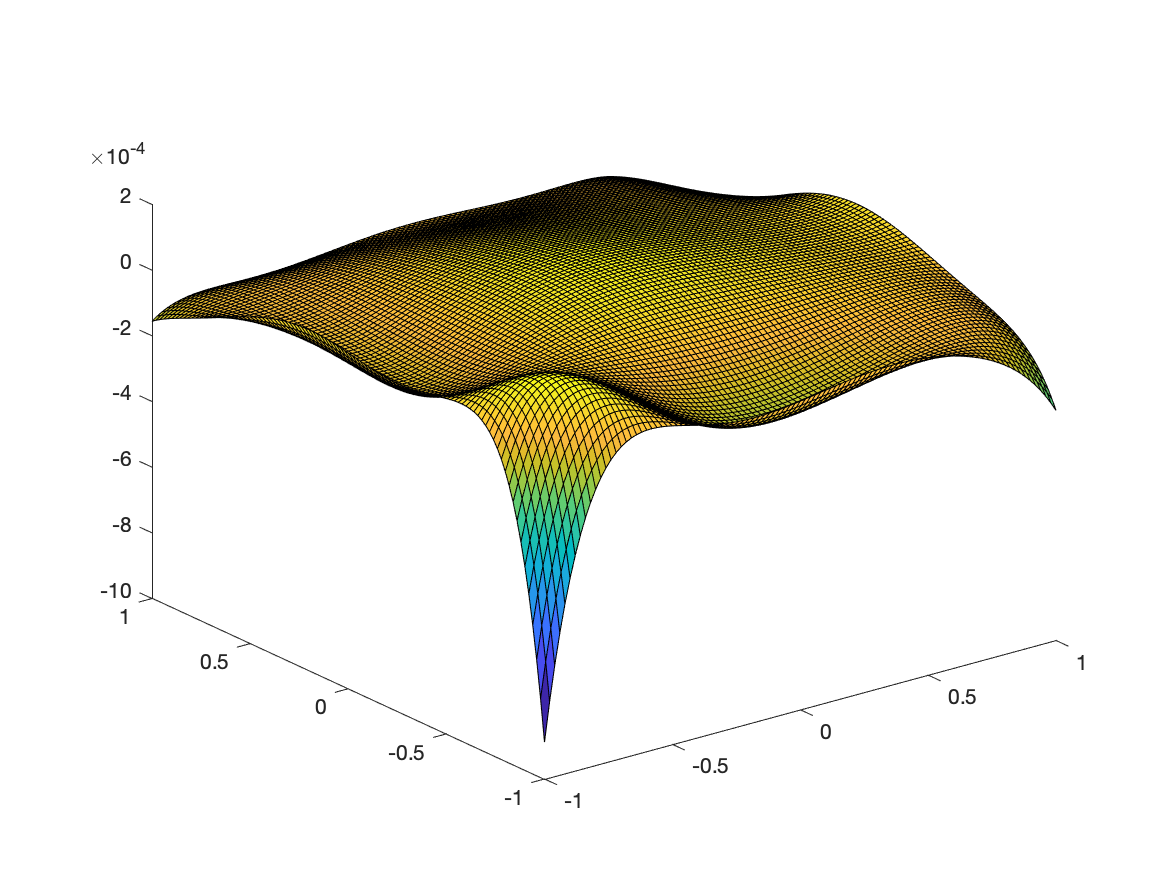} \\ b}
		\end{minipage}
		\hfill
		\begin{minipage}[h]{0.5\linewidth}
			\center{\includegraphics[width=1\linewidth]{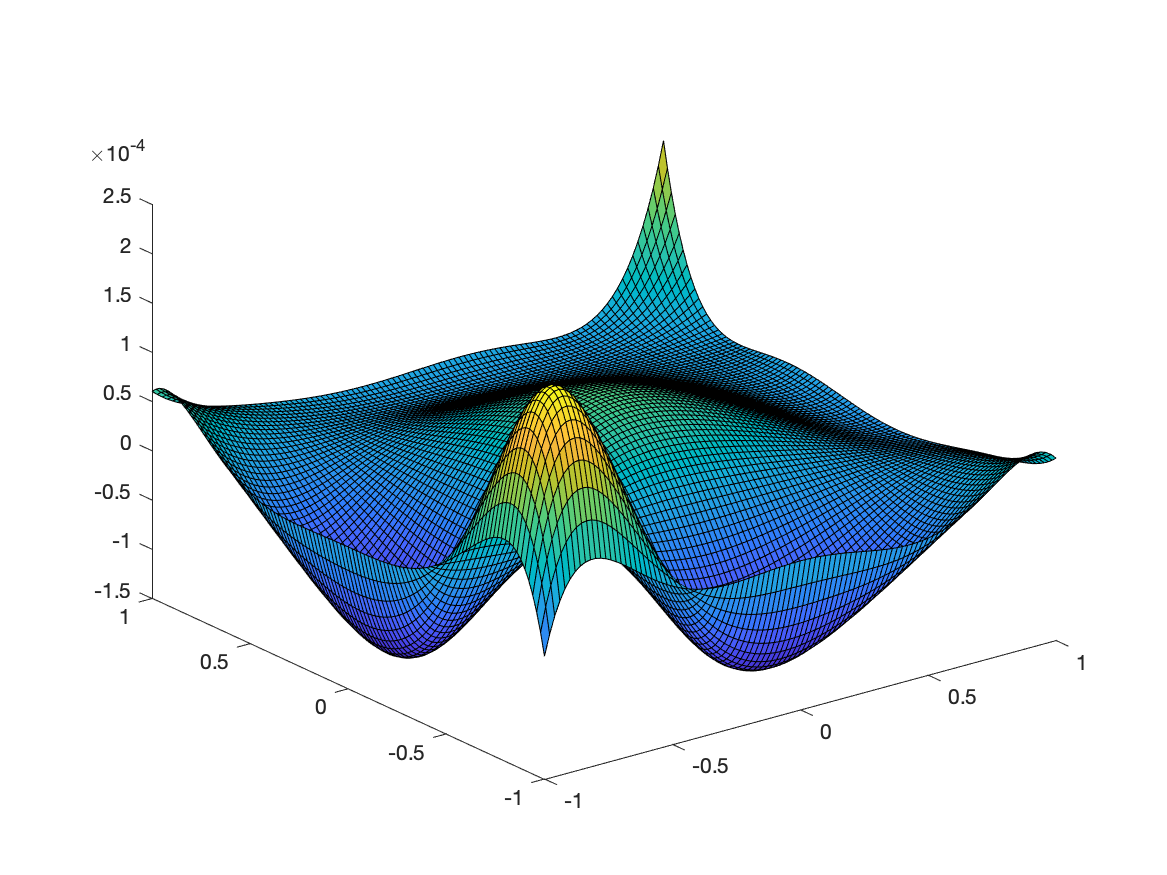} \\ c}
		\end{minipage}
		\begin{minipage}[h]{0.5\linewidth}
			\center{\includegraphics[width=1\linewidth]{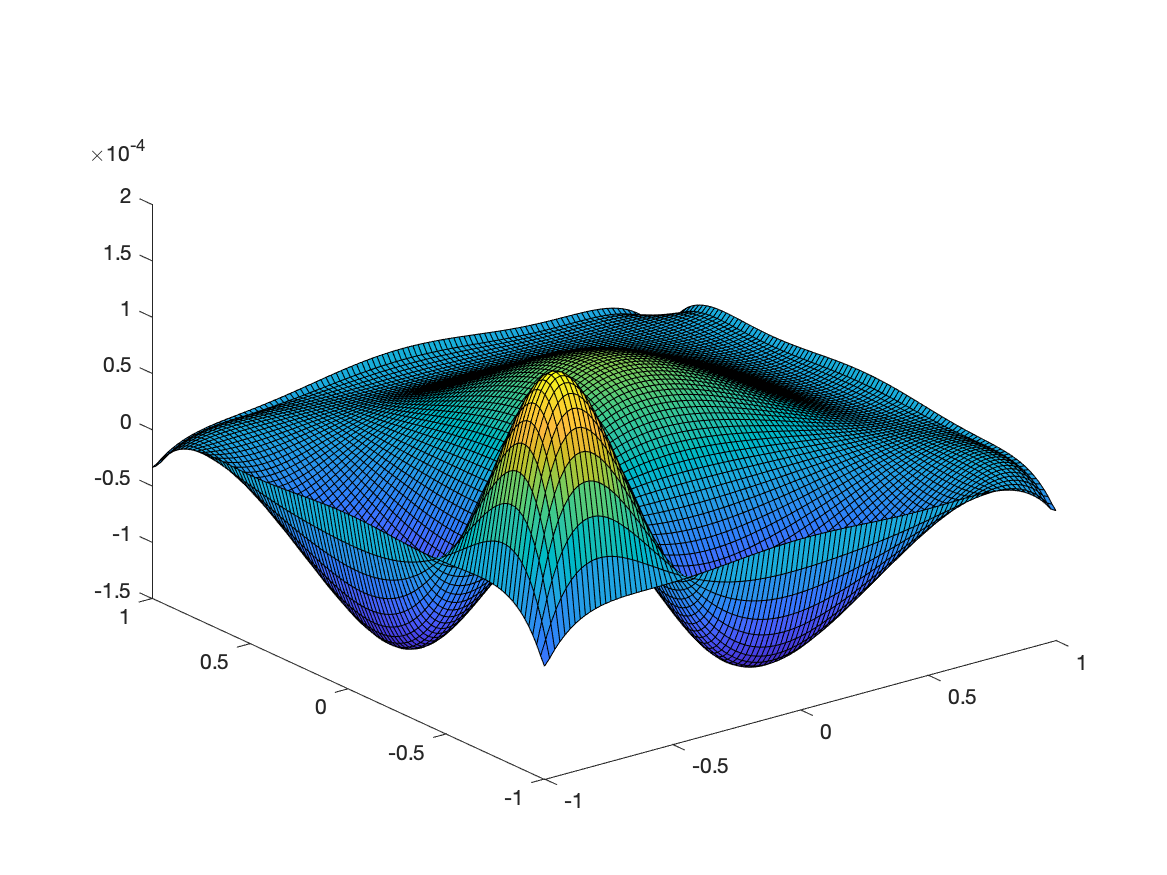} \\ d}
		\end{minipage}
		\caption{Recovery of the derivative $F_1^{(2,2)}$  with  random noise in the  input data . The exact derivative $F_1^{(2,2)}$ (Fig. a ); approximation to $F_1^{(2,2)}$ for    $\delta= 10^{-6}$ (Fig.  b);    for $\delta= 10^{-7}$ (Fig. c) and $\delta= 10^{-8}$ (Fig.  d), }
		\label{Fig2}
	\end{figure}
	\begin{figure}[h!]
		\begin{minipage}[h]{0.5\linewidth}
			\center{\includegraphics[width=1\linewidth]{ExactSpl.eps} \\ a}
		\end{minipage}
		\begin{minipage}[h]{0.5\linewidth}
			\center{\includegraphics[width=1\linewidth]{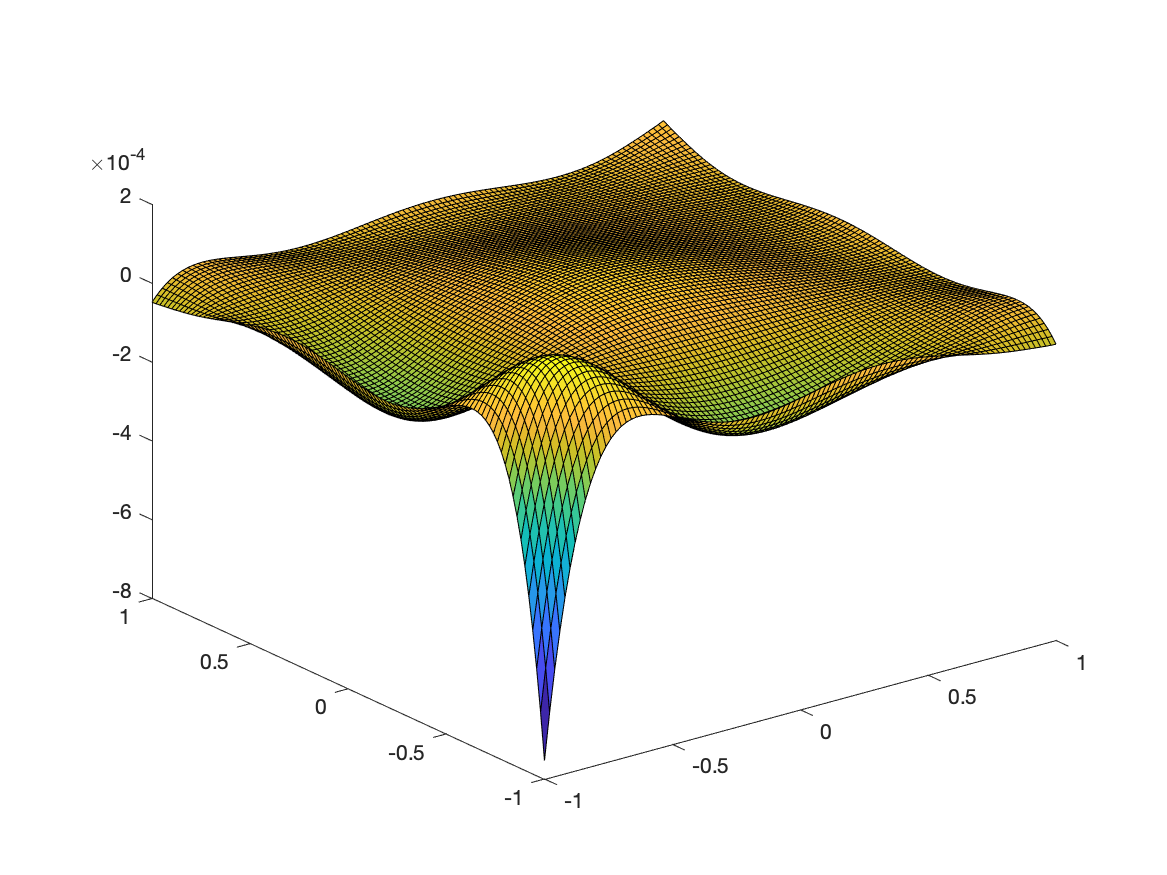} \\ b}
		\end{minipage}
		\hfill
		\begin{minipage}[h]{0.5\linewidth}
			\center{\includegraphics[width=1\linewidth]{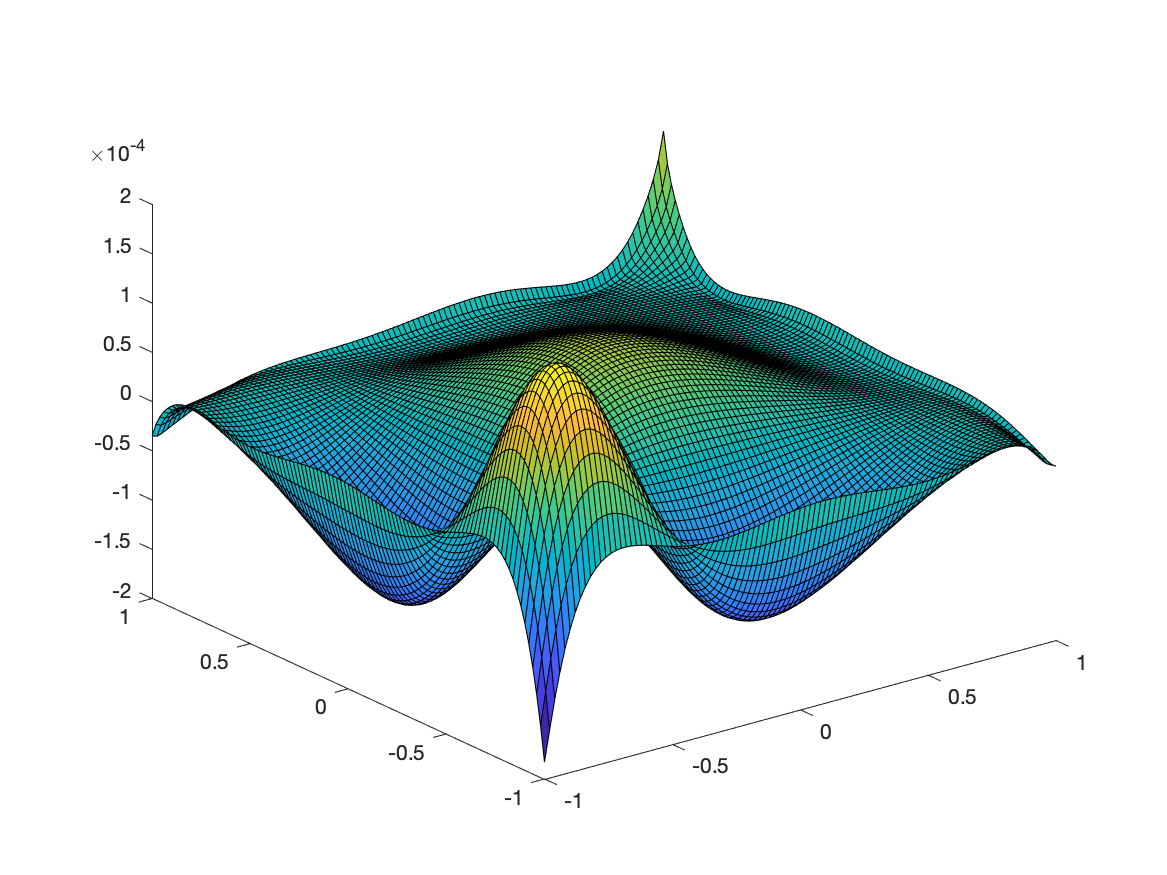} \\ c}
		\end{minipage}
		\begin{minipage}[h]{0.5\linewidth}
			\center{\includegraphics[width=1\linewidth]{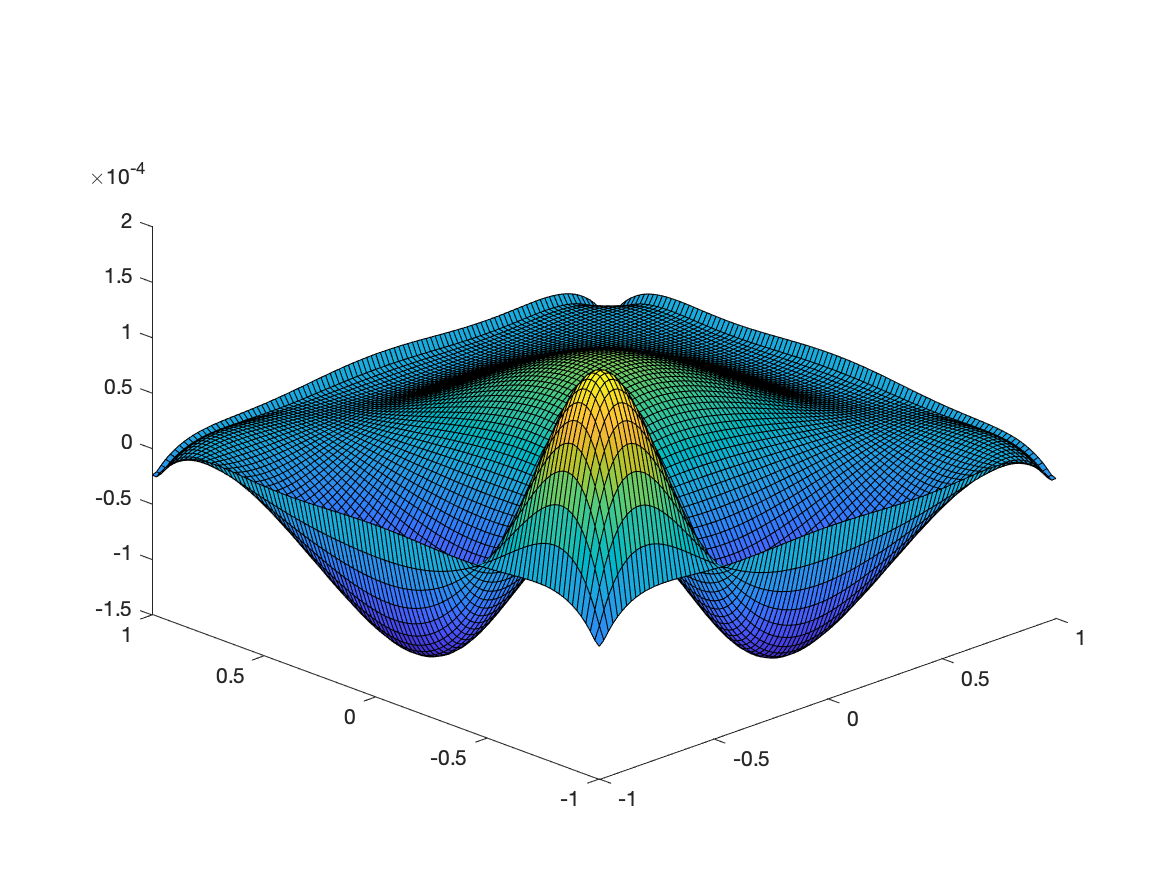} \\ d}
		\end{minipage}
		\caption{Recovery of the derivative $F_1^{(2,2)}$ with noise in the input data arising from the quadrature formula. The exact derivative $F_1^{(2,2)}$ (Fig. a); approximation to $F_1^{(2,2)}$ for    $\delta= 10^{-6}$ (Fig.  b);    for $\delta= 10^{-7}$ (Fig. c) and $\delta= 10^{-8}$ (Fig.  d). }
		\label{Fig1}
	\end{figure}
	As can be seen from the graphs and tables above, for both types of noise, the truncation method gives the same order of accuracy for recovering the derivative $F_1^{(2,2)}$.
	At the same time, applying the quadrature formula expands the area of using the proposed method in computational problems, especially in the situation when the input data are given in the form of a set of function values at the grid nodes.
	
	\subsection{Example 2}
	
	Let us test the method (\ref{ModVer}) on an analytic function. Following \cite{WW2005}, 
	we take the function $ F_2(t,\tau)=(2-(2t-1)^2)^2 \cos(4\tau)/C$.
	Let us put $\mu=6$. It is easy to check that $\|F_2\|_{2,6} \approx 1$ and $\|F_2^{(2,2)}\|_{L_2} \approx 10^{-4}$, if $C=43940129$.
	The Fourier-Legendre coefficients of the considered function are calculated using the quadrature trapezoid formula for $h = 4\cdot 10^{-4}, 10^{-4}, 4\cdot10^{-5}$, which in turn according to formula (\ref{perturbation}) matches $\delta \approx 10^{-6}, 10^{-7}, 10^{-8}$, respectively.

	\begin{table}[h!]
		\centering
		\caption{ The results of recovering derivative $F_2^{(2,2)}$ for noise from quadrature formula }
		\label{tbl3}
		\begin{tabular}{|c|c|c|c|}
			\hline
			$\delta$ & $10^{-6}$ & $10^{-7}$ & $10^{-8}$  \\ \hline
			Error in $L_2$         &  $3.8 \cdot 10^{-5}  $       &  $ 1 \cdot 10^{-6}   $  &  $1.53  \cdot 10^{-7}   $  \\ \hline
			Error in $C$          &  $1,85 \cdot 10^{-4} $       &  $6.37 \cdot 10^{-6}   $  &  $8.17  \cdot 10^{-7}   $  \\ \hline
			n  &  $11$      & $18$    & $25$    \\ \hline
			h  &  $4\cdot 10^{-4}$    & $10^{-4} $  & $ 4\cdot10^{-5}$    \\ \hline
		\end{tabular}
	\end{table}
	
	The results of the numerical experiment are shown in  Table \ref{tbl3} and Graph \ref{Fig3}.
	
	\begin{figure}[h!]
		\begin{minipage}[h]{0.5\linewidth}
			\center{\includegraphics[width=1\linewidth]{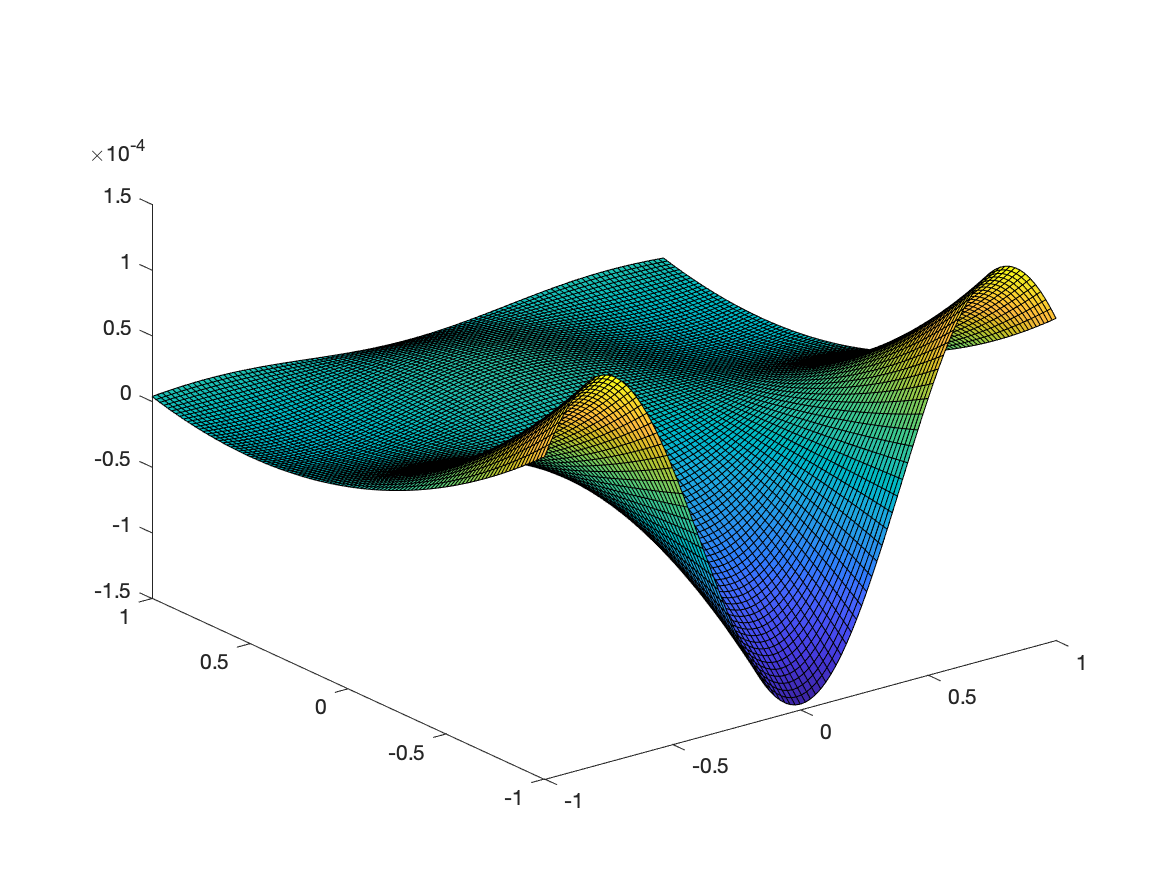} \\ a}
		\end{minipage}
		\begin{minipage}[h]{0.5\linewidth}
			\center{\includegraphics[width=1\linewidth]{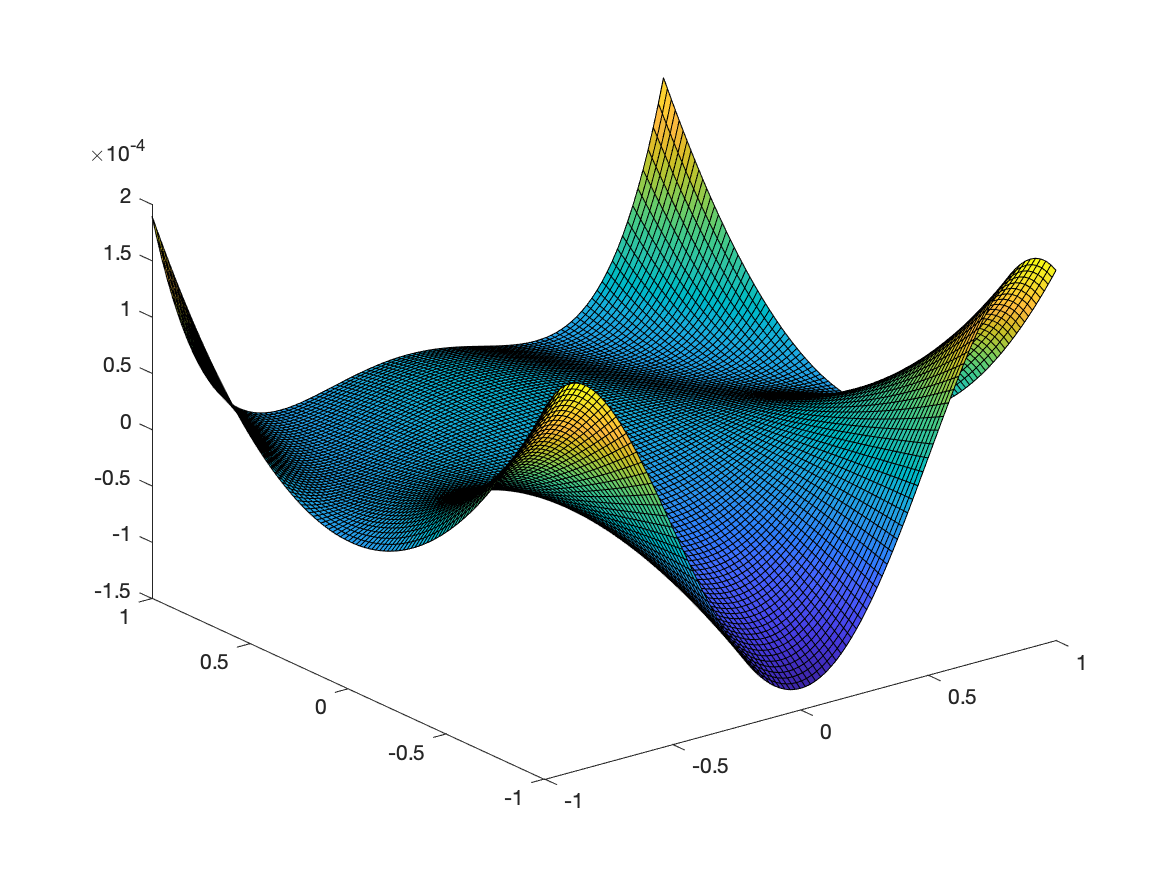} \\ b}
		\end{minipage}
		\hfill
		\begin{minipage}[h]{0.5\linewidth}
			\center{\includegraphics[width=1\linewidth]{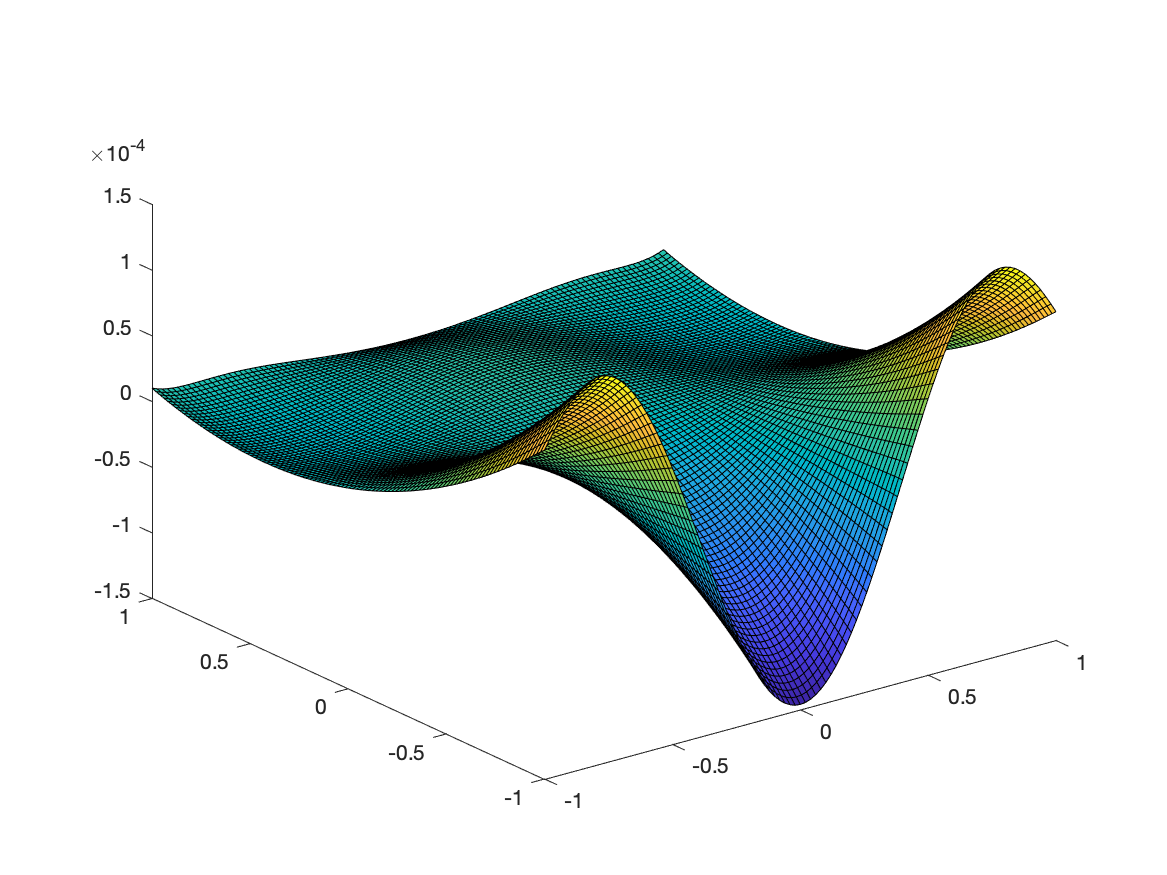} \\ c}
		\end{minipage}
		\begin{minipage}[h]{0.5\linewidth}
			\center{\includegraphics[width=1\linewidth]{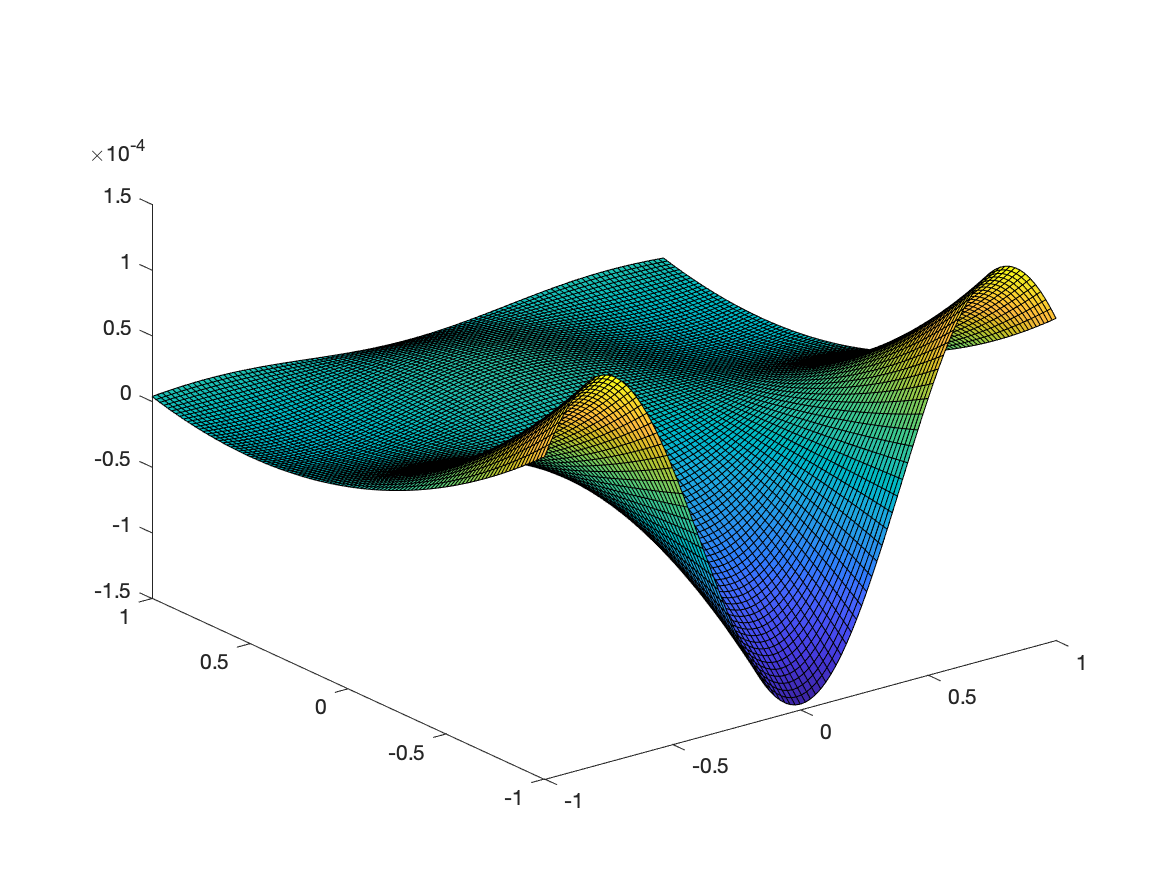} \\ d}
		\end{minipage}
		\caption{Recovery of the derivative $F_2^{(2,2)}$ with noise in the input data arising from the quadrature formula. The exact derivative $F_2^{(2,2)}$ (Fig. a); approximation to $F_2^{(2,2)}$ for    $\delta= 10^{-6}$ (Fig.  b);    for $\delta= 10^{-7}$ (Fig. c) and $\delta= 10^{-8}$ (Fig.  d). }
		\label{Fig3}
	\end{figure}

The numerical results given above show that the proposed method (\ref{ModVer}) works efficiently and coincides well with the theoretical results.

 \section{ Acknowledgements }
 
This project has received funding through the MSCA4Ukraine project,
which is funded by the European Union.
In addition, the first named author is supported by the Volkswagen Foundation project  "From Modeling and Analysis to Approximation".


\begin{thebibliography}{9}





  \bibitem{Ahn&Choi&Ramm_2006}
S. Ahn, U.J. Choi, and  A.G. Ramm, A scheme for stable numerical differentiation,
\emph{J. Comput. Appl. Math.} \textbf{186} (2006),  325--334.


\bibitem{And84} R.S. Anderssen and F.R. de Hoog, Finite difference
methods for the numerical differentiation of non-exact data, \emph{Computing}  \textbf{33} (1984), 259--267.


\bibitem{Cul71} J. Cullum, Numerical Differentiation and Regularization, \emph{SIAM Journal on Numerical Analysis} \textbf{8} (1971), 254--265.

\bibitem{Dolgopolova&Ivanov_USSR_Comput_Math_Math_Phys_1966_Eng}
T.F.~Dolgopolova and V.K.~Ivanov,
On numerical differentiation, \emph{Zh. Vychisl. Mat. and Mat. Ph.}  \textbf{6} (1966),  223--232.

\bibitem{EgorKond_1989}
Y.V. Egorov and V.A. Kondrat'ev, On a problem of numerical differentiation, \emph{
Vestnik Moskov. Univ. Ser. I Mat. Mekh.} \textbf{3} (1989),  80-81.

\bibitem{ErbSem2015}  W. Erb and E.V. Semenova,
On adaptive discretization schemes for the solution of ill-posed problems with semiiterative methods, \emph{Applicable
Analysis}  \textbf{94} (2015), 2057--2076.

\bibitem{Groetsch_1992_V74_N2}
C.W. Groetsch,
Optimal order of accuracy in Vasin's method for differentiation of noisy functions, \emph{J. Optim.Theory Appl.} \textbf{74} (1992), 373--378.

\bibitem{Hanke&Scherzer_2001_V108_N6}
M.~Hanke and  O.~Scherzer,
Inverse problems light: numerical differentiation,
\emph{Amer. Math. Monthly} \textbf{108} (2001),  512--521.

\bibitem{Lu&Naum&Per}
S. Lu, V. Naumova and S.V. Pereverzev,
Legendre polynomials as a recommended basis for numerical differentiation in the
presence of stochastic white noise,
\emph{J. Inverse Ill-Posed Probl.} \textbf{21} (2013), 193--216.

\bibitem{LYHB}
J. Luo, K. Ying and  P. He, J. Bai,
Properties of Savitzky-Golay digital differentiators, \emph{Digital Signal Processing.}
\textbf{15} (2005), 122--136. 


\bibitem{MatPer2002}
P.~Mathe, S.V.~Pereverzev
Stable summation of othogonal series with noisy coefficients, \emph{J. Approx. Theory}  \textbf{118} (2002), 66--80.


\bibitem{Meng&Zhaoa&Mei&Zhou_2020}
Z.~Meng, Z.~Zhaoa, D.~Mei and Y.~Zhou,
Numerical differentiation for two-dimensional functions by a Fourier extension method,
\emph{Inverse Problems in Science and Engineering} \textbf{28} (2020),  1--18.



\bibitem{Mul69}
C. M\" uller, 
\emph{Foundations of the Mathematical Theory of Electromagnetic Waves},
Springer-Verlag, Berlin, 1969.



\bibitem{Mileiko_Solodkii_2017_UMJ}
G.L.~Myleiko and S.G.~Solodky,
Hyperbolic cross and complexity of different classes of linear ill-posed problems, \emph{Ukr. Math. J.} \textbf{69} (2017), 951--963.

\bibitem{Mileiko_Solodkii_2016_ApAn}
G.L.~Myleiko and  S.G.~Solodky, On optimization of projection methods for solving some classes of severely ill-posed
problems, \emph{Applicable Analysis} \textbf{95} (2016), 826--841.

\bibitem{Mileiko_Solodkii_2014}
G.L.~Myleiko and  S.G.~Solodky, The minimal radius of Galerkin information for severely ill-posed problems,
\emph{J. Inverse Ill-Posed Probl.} \textbf{22} (2014),  739--757.





\bibitem{Nakamura&Wang&Wang_2008}
G. Nakamura, S.\, Z. Wang and Y.B. Wang, Numerical differentiation for the second order derivatives of functions of two
variables, \emph{J. Comput. Appl. Math.} \textbf{212} (2008),  341--358.


\bibitem{Pereverzev_Computing_1995}  
S.V. Pereverzev,
Optimization of projection methods for solving ill-posed problems, \emph{Computing} \textbf{55} (1995), 113--124.


\bibitem{PS1996}
S.V. Pereverzev and S.G. Solodky, The minimal radius of Galerkin information for the Fredholm problem of the first
kind, \emph{J. Complexity} \textbf{12} (1996), 401--415.

\bibitem{Qian&Fu&Xiong&Wei_2006}
~Z. Qian,
Fourier truncation method for high order numerical derivatives, \emph{Appl.
Math. Comput.} \textbf{181} (2006), 940--948.


\bibitem{Ramm_1968_No11}
A.G. Ramm,
On numerical differentiation, \emph{Izv. Vuzov. Matem.} \textbf{11} (1968), 131--134.

\bibitem{RammSmir_2001}
A.G. Ramm  and A.B. Smirnova,
On stable numerical differentiation, \emph{Math. Comput.} \textbf{70} (2001),  1131--1153.

\bibitem{Qu96} Q. Ruibin,
A new approach to numerical differentiation and integration, \emph{Mathematical and Computer Modelling}, \textbf{24} (1996) 55--68.

\bibitem{Sem_Sol_2021}
Y.V. Semenova and S.G. Solodky, Error bounds for Fourier-Legendre truncation method in numerical differentiation, \emph{Journal of Numerical and Applied Mathematics} \textbf{3} (2021), 113--130.

\bibitem{Lyashko2022}
Y.V. Semenova and S.G. Solodky,  Optimal methods for recovering mixed derivatives of non-periodic functions. \emph{Journal of Numerical and Applied Mathematics} \textbf{2} (2022), 143--150.

\bibitem{SSS_CMAM}
E.V.Semenova, S.G. Solodky  and S.A. Stasyuk, Application of Fourier Truncation Method to Numerical
Differentiation for Bivariate Functions, \emph{Computational Methods in Applied Mathematics} \textbf{22} (2022), 
477--491. 

\bibitem{SSS_Rew2021} E.V. Semenova, S.G. Solodky  and S.A. Stasyuk,  Truncation method for numerical differentiation problem, \emph{Proceedings of the Institute of Mathematics of the National Academy of Sciences of Ukraine:
	Modern problems of mathematics and its applications} \textbf{18} (2021), 644--672.

\bibitem{Sh2016}
C. Shekhar,
On simplified application of multidimensional Savitzky-Golay filters and differentiators,
\emph{AIP Conference Proceedings}  \textbf{1705} (2015).


\bibitem{SolSha2015}S. G. Solodky and K.K. Sharipov, Summation of smooth functions of two variables with perturbed Fourier coefficients, \emph{J. Inverse Ill-Posed Probl.} \textbf{23}  (2015), 287--297.

\bibitem{Sol_Stas_JC2020} S.G. Solodky and S.A. Stasyuk, Estimates of efficiency for two methods of stable numerical summation
of smooth functions, \emph{J. Complexity} \textbf{56} (2020).






\bibitem{Sol_Stas_UMZ2022}
S.G.  Solodky and S.A. Stasyuk,
Optimization of the Methods of Numerical Differentiation for Bivariate Functions,
\emph{ Ukr. Math. J.} \textbf{74} (2022), 289--313.
%

\bibitem{TrWW} J.F. Traub, G.W. Wasilkowski and H. Wozniakowski, \emph{Information-Based Complexity}, Academic
Press, New York, 1988.

\bibitem{TrW} J.\,F. Traub and H. Wozniakowski, \emph{A General Theory of Optimal Algorithms}, Academic Press, New York,
1980.

\bibitem{VasinVV_1969_V7_N2}
V.V. Vasin, 
Regularization of the numerical differentiation problem
\emph{Ural. Gos. Univ. Mat. Zap.}  \textbf {7} (1969),  29--33.

\bibitem{Wang_Hon_Ch_2006} Y.B. Wang, Y.C. Hon  and  J. Cheng,
Reconstruction of high order derivatives from input data,
\emph{J. Inverse Ill-Posed Probl.} \textbf{14} (2006), 205-218

\bibitem{WW2005}  Y.B.  Wang and T. Wei,
Numerical differentiation for two-dimensional scattered data, \emph{J.Math.Anal. Appl.} \textbf{312} (2005), 121-137.


\bibitem{Zhao_2010}
Z. Zhao,
A truncated Legendre spectral method for solving numerical differentiation
\emph{International Journal of
Computer Mathematics} \textbf{87} (2010), 3209--3217.

\bibitem{Zhao&Meng&Zhao&You&Xie_2016}
Z. Zhao, Z. Meng, L. Zhao, L. You and O. Xie,
A stabilized algorithm for multi-dimensional numerical differentiation
\emph{Journal of Algorithms and Computational Technology} \textbf{10} (2016), 73--81.



\end{thebibliography}

\end{document}